\documentstyle[12pt]{amsart}

\newsymbol\boxtimes 1202

\newcommand{\nc}{\newcommand}

\nc{\ad}{\operatorname{ad}}
\nc{\Aut}{\operatorname{Aut}}
\nc{\Boxtimes}{\fbox{$\times$}}
\nc{\blt}{\bullet}
\nc{\bSt}{\mbox{\bf{St}}}
\nc{\card}{\operatorname{card}}
\nc{\Cch}{\check{C}}
\nc{\cd}{\operatorname{cd}}
\nc{\Ch}{\operatorname{Ch}}
\nc{\chara}{\operatorname{char}}
\nc{\CHom}{\cal{H}om}
\nc{\Coker}{\operatorname{Coker}}
\nc{\codim}{\operatorname{codim}}
\nc{\Cone}{\operatorname{Cone}}
\nc{\cSgn}{\cal{S}gn}
\nc{\depth}{\operatorname{depth}}
\nc{\dirlim}{\underset{\rightarrow}{\operatorname{lim}}}
\nc{\dotbox}{\overset{\bullet}{\boxtimes}}
\nc{\dotimes}{\overset{\bullet}{\otimes}}
\nc{\emp}{\emptyset}
\nc{\Ext}{\operatorname{Ext}}
\nc{\Fac}{\cal{F}ac}
\nc{\Fun}{\operatorname{F}}
\nc{\FS}{\cal{FS}}
\nc{\Hom}{\operatorname{Hom}}
\nc{\hgt}{\operatorname{ht}}
\nc{\Id}{\operatorname{Id}}
\nc{\id}{\operatorname{id}}
\nc{\Ima}{\operatorname{Im}}
\nc{\Ind}{\operatorname{Ind}}
\nc{\infh}{\frac{\infty}{2}}
\nc{\invlim}{\underset{\leftarrow}{\operatorname{lim}}}
\nc{\Ker}{\operatorname{Ker}}
\nc{\Locsys}{\cal{L}ocsys}
\nc{\Ob}{\operatorname{Ob}}
\nc{\Or}{\cal{O}r}
\nc{\Ord}{\cal{O}rd}
\nc{\Part}{\cal{P}art}
\nc{\PGL}{\operatorname{PGL}}
\nc{\sgn}{\operatorname{sgn}}
\nc{\Sh}{\cal{S}h}
\nc{\supp}{\operatorname{supp}}
\nc{\tFS}{\widetilde{\cal{FS}}}
\nc{\Tor}{\operatorname{Tor}}
\nc{\totimes}{\tilde{\otimes}}
\nc{\Vect}{\cal{V}ect}
\nc{\wt}{\widetilde}

\nc{\tDFS}{{\widetilde{\cal{DFS}}}}
\nc{\DFS}{{\cal{DFS}}}
\nc{\bHom}{\overline{\operatorname{Hom}}}
\nc{\DVect}{\cal{DV}ect}

\nc{\bo}{\mbox{\bf{0}}}
\nc{\One}{\mbox{\bf{1}}}
\nc{\one}{\mbox{\bf{1}}}

\nc{\BA}{\Bbb A}
\nc{\ba}{\mbox{\bf{a}}}
\nc{\baJ}{\bar{J}}
\nc{\BAO}{\overset{\circ}{\Bbb A}}
\nc{\BB}{\Bbb B}
\nc{\bB}{\mbox{\bf{B}}}
\nc{\BC}{\Bbb C}
\nc{\bCC}{\bar{\cal{C}}}
\nc{\bD}{\bar{D}}
\nc{\bd}{\mbox{\bf{d}}}
\nc{\BE}{\overline{E}}
\nc{\BF}{\overline{F}}
\nc{\bF}{\mbox{\bf{F}}}
\nc{\bL}{\mbox{\bf{L}}}
\nc{\blambda}{\bar{\lambda}}
\nc{\bM}{\mbox{\bf{M}}}
\nc{\bmu}{\vec{\mu}}
\nc{\BN}{\Bbb N}
\nc{\bnu}{\vec{\nu}}
\nc{\BP}{\Bbb P}
\nc{\bP}{\mbox{\bf{P}}}
\nc{\BPO}{\overset{\circ}{\Bbb{P}}}
\nc{\BQ}{\Bbb Q}
\nc{\bq}{\mbox{\bf{q}}}
\nc{\BR}{\Bbb R}
\nc{\br}{\mbox{\bf{r}}}
\nc{\breta}{\bar{\eta}}
\nc{\bs}{\mbox{\bf{s}}}
\nc{\bt}{\mbox{\bf{t}}}
\nc{\bU}{\mbox{\bf{U}}}
\nc{\bu}{\mbox{\bf{u}}}
\nc{\BUpsilon}{\bar{\Upsilon}}
\nc{\bw}{\mbox{\bf{w}}}
\nc{\bx}{\mbox{\bf{x}}}
\nc{\BZ}{\Bbb Z}
\nc{\bz}{\mbox{\bf{z}}}

\nc{\CA}{\cal A}
\nc{\CAD}{\overset{\bullet}{\cal{A}}}
\nc{\CAO}{\overset{\circ}{\cal{A}}}
\nc{\CB}{\cal B}
\nc{\CC}{\cal C}
\nc{\CD}{\cal D}
\nc{\CE}{\cal E}
\nc{\CF}{\cal F}
\nc{\CH}{\cal H}
\nc{\CI}{\cal I}
\nc{\CID}{\overset{\bullet}{\cal{I}}}
\nc{\CJ}{\cal J}
\nc{\CK}{\cal K}
\nc{\CL}{\cal L}
\nc{\CM}{\cal M}
\nc{\CN}{\cal N}
\nc{\CO}{\cal O}
\nc{\CP}{\cal P}
\nc{\CPO}{\overset{\circ}{\cal{P}}}
\nc{\CQ}{\cal Q}
\nc{\CR}{\cal R}
\nc{\CS}{\cal S}
\nc{\CT}{\cal T}
\nc{\CTD}{\overset{\bullet}{\cal{T}}}
\nc{\CTPO}{\overset{\circ}{\cal{T}\cal{P}}}
\nc{\CU}{\cal{U}}
\nc{\CV}{\cal V}
\nc{\CX}{\cal X}
\nc{\CY}{\cal Y}
\nc{\CZ}{\cal Z}

\nc{\DO}{\overset{\circ}{D}}
\nc{\dpar}{\partial}
\nc{\DT}{\overset{\bullet}{T}}

\nc{\fA}{\frak{A}}
\nc{\fE}{\frak{E}}
\nc{\fF}{\frak F}
\nc{\ff}{\frak f}
\nc{\fg}{\frak g}
\nc{\fl}{\frak{l}}
\nc{\fp}{\frak{p}}
\nc{\fu}{\frak{u}}

\nc{\HO}{\overset{\circ}{H}}
\nc{\hfg}{\hat{\frak{g}}}
\nc{\hL}{\hat{L}}

\nc{\jo}{\overset{\circ}{j}}

\nc{\phid}{\overset{\bullet}{\phi}}

\nc{\tBP}{\tilde{\Bbb{P}}}
\nc{\tC}{\tilde{C}}
\nc{\tc}{\tilde{c}}
\nc{\tCA}{\tilde{\cal{A}}}
\nc{\tCC}{\tilde{\cal{C}}}
\nc{\tCI}{\tilde{\cal{I}}}
\nc{\tCO}{\tilde{\cal{O}}}
\nc{\tCP}{\tilde{\cal{P}}}
\nc{\tCT}{\tilde{\cal{T}}}
\nc{\tD}{\tilde{D}}
\nc{\tDelta}{\tilde{\Delta}}
\nc{\tE}{\tilde E}
\nc{\tF}{\tilde F}
\nc{\tfF}{\tilde{\frak{F}}}
\nc{\tfu}{\tilde{\frak{u}}}
\nc{\tJ}{\tilde{J}}
\nc{\tj}{\tilde{j}}
\nc{\tK}{\tilde K}
\nc{\tL}{\tilde{L}}
\nc{\tM}{\tilde{M}}
\nc{\tP}{\tilde{P}}
\nc{\tPhi}{\tilde{\Phi}}
\nc{\TPO}{\overset{\circ}{T\Bbb{P}}}
\nc{\tR}{\tilde{R}}
\nc{\tS}{\tilde S}
\nc{\tT}{\tilde{T}}
\nc{\ttau}{\tilde{\tau}}
\nc{\ttheta}{\tilde{\theta}}
\nc{\tU}{\tilde{U}}
\nc{\tUpsilon}{\tilde{\Upsilon}}
\nc{\ty}{\tilde y}
\nc{\tY}{\tilde Y}
\nc{\txi}{\tilde{\xi}}

\nc{\UD}{\overset{\bullet}{U}}
\nc{\UO}{\overset{\circ}{U}}

\nc{\valpha}{\vec{\alpha}}
\nc{\vbeta}{\vec{\beta}}
\nc{\vc}{\vec{c}}
\nc{\vD}{\vec{D}}
\nc{\vd}{\vec{d}}
\nc{\vgamma}{\vec{\gamma}}
\nc{\vlambda}{\vec{\lambda}}
\nc{\vmu}{\vec{\mu}}
\nc{\vnu}{\vec{\nu}}
\nc{\vo}{\vec{0}}
\nc{\vu}{\vec{u}}
\nc{\vx}{\vec{x}}

\nc{\XO}{\overset{\circ}{X}}

\nc{\nen}{\newenvironment}
\nc{\ol}{\overline}
\nc{\ul}{\underline}
\nc{\ra}{\rightarrow}
\nc{\lra}{\longrightarrow}
\nc{\Lra}{\Longrightarrow}
\nc{\lla}{\longleftarrow}
\nc{\Llra}{\Longleftrightarrow}
\nc{\hra}{\hookrightarrow}
\nc{\iso}{\overset{\sim}{\lra}}

\nc{\Thm}[1]{Theorem~\ref{#1}}
\nc{\Prop}[1]{Proposition~\ref{#1}}
\nc{\Lem}[1]{Lemma~\ref{#1}}
\nc{\Cor}[1]{Corollary~\ref{#1}}
\nc{\Conj}[1]{Conjecture~\ref{#1}}
\nc{\Claim}[1]{Claim~\ref{#1}}
\nc{\Defn}[1]{Definition~\ref{#1}}
\nc{\Exa}[1]{Example~\ref{#1}}
\nc{\Rem}[1]{Remark~\ref{#1}}
\nc{\Note}[1]{Note~\ref{#1}}

\nen{thm}[1]{\label{#1}{\bf Theorem.\ } \em}{}
\nen{prop}[1]{\label{#1}{\bf Proposition.\ } \em}{}
\nen{lem}[1]{\label{#1}{\bf Lemma.\ } \em}{}
\nen{cor}[1]{\label{#1}{\bf Corollary.\ } \em}{}
\nen{conj}[1]{\label{#1}{\bf Conjecture.\ } \em}{}
\nen{claim}[1]{\label{#1}{\bf Claim.\ } \em}{}

\nen{defn}[1]{\label{#1}{\bf Definition.\ } }{}
\nen{exa}[1]{\label{#1}{\bf Example.\ } }{}

\nen{rem}[1]{\label{#1}{\em Remark.\ } }{}
\nen{note}[1]{\label{#1}{\em Note.\ } }{}
\nen{exer}[1]{\label{#1}{\em Exercise.\ } }{}

\begin{document}

\title[]{Localization of $\fu$-modules. IV.\\
Localization on $\BP^1$}
\author{Michael Finkelberg}
\address{Independent Moscow University, 65-3 Mikloukho-Maklai St.,
apt. 86, Moscow 117342 Russia}
\email{fnklberg@@ium.ips.ras.ru}
\author{Vadim Schechtman}
\address{Dept. of Mathematics, SUNY at Stony Brook, Stony Brook,
NY 11794-3651 USA}
\email{vadik@@math.sunysb.edu}
\thanks{The first author was partially supported by
grants from AMS, SOROS, INTAS, NSF,
and Paul and Gabriela Rosenbaum Foundation.
The second author was partially supported by NSF grant DMS-9202280}
\date{June 1995; revised: September 1995, November 1995, December 1995\\
q-alg/9506011}
\maketitle

\section{Introduction}

\subsection{} This article is a sequel to ~\cite{fs}.
Given a collection of $m$ finite factorizable sheaves $\{\CX_k\}$,
we construct here some perverse sheaves over configuration
spaces of points on a projective line $\BP^1$ with $m$ additional
marked points.

We announce here (with sketch proof) the computation of the
cohomology spaces of these sheaves. They turn out
to coincide with certain "semiinfinite" $\Tor$ spaces of the corresponding
$\fu$-modules introduced by S.Arkhipov. For a precise formulation see
Theorem ~\ref{global thm}.

This result is strikingly similar to the following hoped-for
picture of affine Lie algebra representation theory, explained to us by
A.Beilinson. Let $M_1,M_2$ be two modules over an affine Lie algebra $\hfg$
on the critical level. One hopes that there is a localization functor
which associates to these modules perverse sheaves $\Delta(M_1),
\Delta(M_2)$ over the semiinfinite flag space $\hat{G}/\hat{B}^0$. Suppose that
$\Delta(M_1)$ and $\Delta(M_2)$ are equivariant with respect to the opposite
Borel subgroups of $\hat{G}$. Then the intersection $S$ of their supports
is finite dimensional, and one hopes that
$$
R^{\blt}\Gamma(S,\Delta(M_1)\otimes\Delta(M_2))=
\Tor_{\infh-\blt}^{\hfg}(M_1,M_2)
$$
where in the right hand side we have
the Feigin (Lie algebra) semiinfinite homology.

As a corollary of Theorem ~\ref{global thm} we get a description of
local systems of conformal blocks in WZW models in genus zero
(cf. ~\cite{ms}) as natural subquotients of some semisimple local systems
of geometric origin. In particular, these local systems are
semisimple themselves.

In the next paper of the series (joint with R.Bezroukavnikov) we will
explain how to glue factorizable sheaves into punctured curves
of arbitrary genus. This will give rise to an example of a modular
functor. In fact, the cohesive local system described in Chapter 1
of this paper is a genus $0$ case of the general construction
due to R.Bezroukavnikov.

\subsection{} We are grateful to S.Arkhipov and G.Lusztig for the permission
to use their unpublished results. Namely, Theorem ~\ref{braid thm}
about braiding in the category $\CC$ is due to G.Lusztig, and
Chapter 2 (semiinfinite homological algebra in $\CC$) is an exposition
of the results due to S.Arkhipov.

We are also grateful to A.Kirillov, Jr. who explained to us how
to handle the conformal blocks of non simply laced Lie algebras.

\subsection{}
Unless specified otherwise, we will keep the notations of ~\cite{fs}.
For $\alpha=\sum_i\ a_ii\in\BN[I]$ we will use the notation
$|\alpha|:=\sum_i\ a_i$.

References to {\em loc.cit.} will look like Z.1.1 where Z=I, II or III.

We will keep assumptions of III.1.4 and III.16. In particular,
a "quantization" parameter $\zeta$ will be a primitive $l$-th root
of unity where $l$ is a fixed positive number prime to $2,3$.

\newpage
\begin{center}
{\bf CHAPTER 1. Gluing over $\BP^1$}
\end{center}
\vspace{.8cm}

\section{Cohesive local system}

\subsection{Notations} Let $\alpha\in\BN[X]$,
$\alpha=\sum a_\mu\mu$.
We denote by $\supp\alpha$ the subset of $X$ consisting of all $\mu$ such that
$a_\mu\not=0$.
Let  $\pi:\ J\lra X$ be an {\em unfolding} of $\alpha$, that is
a map of sets such that
$|\pi^{-1}(\mu)|=a_\mu$ for any $\mu\in X$.
As always, $\Sigma_\pi$ denotes the group of
automorphisms of $J$ preserving the fibers of $\pi$.

$\BP^1$ will denote a complex projective line. The $J$-th cartesian power
$\BP^{1J}$ will be denoted by $\CP^J$. The group $\Sigma_\pi$ acts
naturally on $\CP^J$, and the quotient
space $\CP^J/\Sigma_\pi$ will be denoted by $\CP^{\alpha}$.

$\CP^{oJ}$ (resp., $\CP^{o\alpha}$) stands for the complement
to diagonals in $\CP^J$ (resp., in $\CP^{\alpha}$).

$T\CP^J$ stands for the total space of the tangent bundle to $\CP^J$;
its points are couples $((x_j),(\tau_j))$ where $(x_j)\in\CP^J$ and
$\tau_j$ is a tangent vector at $x_j$. An open subspace
$$
\DT\CP^J\subset T\CP^J
$$
consists of couples with $\tau_j\neq 0$ for all $j$.
So, $\DT\CP^J\lra\CP^J$ is a $(\BC^*)^J$-torsor.
We denote by $T\CP^{oJ}$ its restriction to $\CP^{oJ}$.
The group $\Sigma_\pi$ acts freely on $T\CP^{oJ}$, and we denote the quotient
$T\CP^{oJ}/\Sigma_\pi$ by $T\CP^{o\alpha}$.

The natural projection
$$
T\CP^{oJ}\lra T\CP^{o\alpha}
$$
will be denoted by $\pi$, or sometimes by $\pi_J$.


\subsection{}
Let $\BP^1_{st}$ ($st$ for "standard") denote "the"
projective line with fixed coordinate $z$; $D_{\epsilon}\subset\BP^1_{st}$ ---
the open disk of radius $\epsilon$ centered at $z=0$; $D:=D_1$.
We will also use the notation $D_{(\epsilon,1)}$ for the open
annulus $D-\overline{D}_{\epsilon}$ (bar means the closure).

The definitions of $D^J, D^{\alpha}, D^{oJ}, D^{o\alpha}, \DT D^J, TD^{oJ},
TD^{o\alpha}$,
etc., copy the above definitions, with $D$ replacing $\BP^1$.

\subsection{} Given a finite set $K$, let
$\tCP^K$ denote the space of $K$-tuples $(u_k)_{k\in K}$ of algebraic
isomorphisms $\BP^1_{st}\iso\BP^1$ such that the images $u_k(D)$
do not intersect.

Given a $K$-tuple $\valpha=(\alpha_k)\in\BN[X]^K$ such that
$\alpha=\sum_k\ \alpha_k$, define a space
$$
\CP^{\valpha}:=\tCP^K\times\prod_{k\in K}\ \DT D^{\alpha_k}
$$
and an open subspace
$$
\CP^{o\valpha}:=\tCP^K\times\prod_{k\in K}\ TD^{o\alpha_k}.
$$
We have an evident "substitution" map
$$
q_{\valpha}:\CP^{\valpha}\lra\DT\CP^{\alpha}
$$
which restricts to $q_{\valpha}:\CP^{o\valpha}\lra T\CP^{o\alpha}$.

\subsubsection{} In the same way we define the spaces
$TD^{\valpha}$, $TD^{o\valpha}$.

\subsubsection{}
\label{assoc space} Suppose that we have an epimorphism
$\xi:L\lra K$, denote $L_k:=\xi^{-1}(k)$. Assume that
each $\alpha_k$ is in turn decomposed
as
$$
\alpha_k=\sum_{l\in L_k} \ \alpha_{l},\ \alpha_{l}\in\BN[X];
$$
set $\valpha_k=(\alpha_{l})\in\BN[X]^{L_k}$. Set
$\valpha_L=(\alpha_l)\in\BN[X]^L$.

Let us define spaces
$$
\CP^{\valpha_L;\xi}=\tCP^{K}\times\prod_{k\in K}\ \DT D^{\valpha_k}
$$
and
$$
\CP^{o\valpha_L;\xi}=\tCP^{K}\times\prod_{k\in K}\ TD^{o\valpha_k}.
$$
We have canonical substitution maps
$$
q_{\valpha_L;\xi}^1:\CP^{\valpha_L;\xi}\lra\CP^{\valpha}
$$
and
$$
q_{\valpha_L;\xi}^2:\CP^{\valpha_L;\xi}\lra\CP^{\valpha_L}.
$$
Obviously,
$$
q_{\valpha}\circ q_{\valpha_L;\xi}^1=q_{\valpha_L}\circ q_{\valpha_L;\xi}^2.
$$

\subsection{Balance function}
\label{fun n}
Consider a function $n:\ X\lra\BZ[\frac{1}{2\det A}]$ such that $$
n(\mu+\nu)=n(\mu)+n(\nu)+\mu\cdot\nu
$$
It is easy to see that $n$ can be written in the following form:
$$
n(\mu)=\frac{1}{2}\mu\cdot\mu+\mu\cdot\nu_0
$$
for some $\nu_0\in X$.
{}From now on we fix such a function $n$ and hence the corresponding
$\nu_0$.

\subsection{} For an arbitrary $\alpha\in\BN[X]$,
let us define a one-dimensional local system
$\CI^{\alpha}_D$ on $TD^{o\alpha}$. We will proceed in the same way as in
III.3.1.

Pick an unfolding of $\alpha$, $\pi:J\lra X$. Define a local system
$\CI^J_D$ on $TD^{oJ}$ as follows: its stalk at each point
$((\tau_j),(x_j))$ where all $x_j$ are real, and all the tangent vectors
$\tau_j$ are real and directed to the right, is $B$. Monodromies are:

--- $x_i$ moves counterclockwise around $x_j$: monodromy is
    $\zeta^{-2\pi(i)\cdot\pi(j)}$;

--- $\tau_j$ makes a counterclockwise circle: monodromy is
    $\zeta^{-2n(\pi(j))}$.

This local system has an evident $\Sigma_{\pi}$-equivariant structure,
and we define a local system $\CI^{\alpha}_D$ as
$$
\CI^{\alpha}_D=(\pi_*\CI^J_D)^{\sgn}
$$
where $\pi:\ TD^{oJ}\lra TD^{o\alpha}$ is the canonical projection,
and $(\bullet)^{\sgn}$ denotes the subsheaf of skew $\Sigma_{\pi}$-invariants.

\subsection{} We will denote the unique homomorphism
$$
\BN[X]\lra X
$$
identical on $X$, as $\alpha\mapsto\alpha^{\sim}$.

\subsubsection{}
\label{admis}
{\bf Definition.} {\em An element $\alpha\in\BN[X]$
is called {\em admissible} if $\alpha^{\sim}\equiv -2\nu_0\bmod lY$.}

\subsection{} We have a canonical "$1$-jet at $0$" map
$$
p_K:\tCP^K\lra T\CP^{oK}
$$

\subsection{Definition} {\em A {\em cohesive local system} (CLS)
(over $\BP^1$) is the following collection of data:

(i) for each admissible $\alpha\in\BN[X]$ a one-dimensional
local system $\CI^{\alpha}$ over $T\CP^{o\alpha}$;

(ii) for each decomposition $\alpha=\sum_{k\in K}\ \alpha_k$, $\alpha_k\in
\BN[X]$, a
factorization isomorphism
$$
\phi_{\valpha}:q^*_{\valpha}\CI^{\alpha}\iso
p_K^*\pi_K^*\CI^{\alpha_K}\boxtimes\Boxtimes_k\ \CI^{\alpha_k^{\sim}}_D
$$
Here $\alpha_K:=\sum_k\ \alpha_k^{\sim}\in\BN[X]$ (note that
$\alpha_K$ is obviously admissible); $\pi_K: T\CP^{oK}\lra T\CP^{o\alpha}$
is the symmetrization map.

These isomorphisms must satisfy the following

{\em Associativity axiom.} In the assumptions of ~\ref{assoc space} the
equality
$$
\phi_{\valpha_L;\xi}\circ
q_{\valpha_L;\xi}^{1*}(\phi_{\valpha})=q_{\valpha_L;\xi}^{2*}(\phi_{\valpha_L})
$$
should hold. Here $\phi_{\valpha_L;\xi}$ is induced by the evident
factorization isomorphisms for local systems on the disk $\CI^{\valpha_k}_D$.}

Morphisms between CLS's are defined in the obvious way.

\subsection{Theorem.} {\em Cohesive local systems over $\BP^1$
exist. Every two CLS's are isomorphic.
The group of automorphisms of a CLS is $B^*$.}

This theorem is a particular case of a more general theorem, valid
for curves of arbitrary genus, to be proved in Part V. We leave the proof
in the case of $\BP^1$ to the interested reader.

\section{Gluing}

\subsection{}
\label{bal fun} Let us define an element $\rho\in X$ by the condition
$\langle \rho,i\rangle=1$ for all $i\in I$. From now on we choose a balance
function $n$, cf. ~\ref{fun n}, in the form
$$
n(\mu)=\frac{1}{2}\mu\cdot\mu + \mu\cdot\rho.
$$
It has the property that $n(-i')=0$ for all $i\in I$. Thus, in the notations
of {\em loc. cit.} we set
$$
\nu_0=\rho.
$$

We pick a corresponding CLS $\CI=\{\CI^{\beta},\ \beta\in\BN[X]\}$.

Given $\alpha=\sum a_ii\in \BN[I]$ and $\vmu=(\mu_k)\in X^K$, we define an
element
$$
\alpha_{\vmu}=\sum a_i\cdot (-i')+\sum_k\mu_k\in\BN[X]
$$
where the sum in the right hand side is a formal one. We say that a pair
$(\vmu,\alpha)$ is {\em admissible} if $\alpha_{\vmu}$ is admissible in
the sense of the previous section, i.e.
$$
\sum_k\ \mu_k-\alpha\equiv -2\rho\bmod lY.
$$
Note that given $\vmu$, there exists $\alpha\in\BN[I]$ such that
$(\vmu,\alpha)$
is admissible if and only if $\sum_k\ \mu_k\in Y$; if this holds true,
such elements $\alpha$ form an obvious countable set.

We will denote by
$$
e:\BN[I]\lra \BN[X]
$$
a unique homomorphism sending $i\in I$ to $-i'\in X$.

\subsection{} Let us consider the space $\DT\CP^K\times\DT\CP^{e(\alpha)}$;
its points are quadruples $((z_k),(\tau_k),(x_j),(\omega_j))$ where
$(z_k)\in\CP^K,\ \tau_k$ --- a non-zero tangent vector to $\BP^1$ at
$z_k$, $(x_j)\in\CP^{e(\alpha)}$, $\omega_j$ --- a non-zero tangent
vector at $x_j$. To a point $z_k$ is assigned a weight $\mu_k$,
and to $x_j$ --- a weight $-\pi(j)'$. Here $\pi:J\lra I$ is an unfolding
of $\alpha$ (implicit in the notation $(x_j)=(x_j)_{j\in J}$).

We will be interested in some open subspaces:
$$
\DT\CP^{\alpha}_{\vmu}:=T\CP^{oK}\times\DT\CP^{e(\alpha)}\subset
\DT\CP^{\alpha_{\vmu}}
$$
and
$$
T\CP^{o\alpha}_{\vmu}\subset
\DT\CP^{\alpha}_{\vmu}
$$
whose points are quadruples $((z_k),(\tau_k),(x_j),(\omega_j))\in
\DT\CP^{\alpha}_{\vmu}$  with all $z_k\neq x_j$. We have an obvious
symmetrization projection
$$
p^{\alpha}_{\vmu}: T\CP^{o\alpha}_{\vmu}\lra T\CP^{o\alpha_{\vmu}}.
$$
Define a space
$$
\CP^{\alpha}_{\vmu}=T\CP^{oK}\times\CP^{e(\alpha)};
$$
its points are triples $((z_k),(\tau_k),(x_j))$ where
$(z_k)$, $(\tau_k)$ and $(x_j)$ are as above; and to $z_k$ and $x_j$
the weights as above are assigned. We have the canonical projection
$$
\DT\CP^{\alpha}_{\vmu}\lra \CP^{\alpha}_{\vmu}.
$$
We define the open subspaces
$$
\CP^{o\alpha}_{\vmu}\subset
\CP^{\bullet\alpha}_{\vmu}\subset
\CP^{\alpha}_{\vmu}.
$$
Here the $\bullet$-subspace (resp., $o$-subspace)
consists of all $((z_k),(\tau_k),(x_j))$ with
$z_k\neq x_j$ for all $k,j$ (resp., with all $z_k$ and $x_j$ distinct).

We define the {\em principal stratification} $\CS$ of $\CP^{\alpha}_{\vmu}$
as the stratification generated by subspaces $z_k=x_j$ and $x_j=x_{j'}$
with $\pi(j)\neq\pi(j')$. Thus, $\CP^{o\alpha}_{\vmu}$ is the open
stratum of $\CS$. As usually, we will denote by the same letter the
induced stratifications on subspaces.

The above projection restricts to
$$
T\CP^{o\alpha}_{\vmu}\lra \CP^{o\alpha}_{\vmu}.
$$

\subsection{Factorization structure}
\subsubsection{} Suppose we are given $\valpha\in\BN[I]^K,\ \beta\in\BN[I]$;
set $\alpha:=\sum_k\ \alpha_k$. Define a space
$$
\CP^{\valpha,\beta}_{\vmu}\subset
\tCP^K\times\prod_k\ D^{\alpha_k}\times\CP^{e(\beta)}
$$
consisting of all collections $((u_k),((x^{(k)}_j)_k),(y_j))$
where $(u_k)\in\tCP^K,\ (x^{(k)}_j)_k\in D^{\alpha_k},\
(y_j)\in\CP^{e(\beta)}$, such that
$$
y_j\in\BP^1-\bigcup_{k\in K}\ \overline{u_k(D)}
$$
for all $j$ (the bar means closure).

We have canonical maps
$$
q_{\valpha,\beta}:\CP^{\valpha,\beta}_{\vmu}\lra\CP^{\alpha+\beta}_{\vmu},
$$
assigning to $((u_k),((x^{(k)}_j)_k),(y_j))$ a configuration
$(u_k(0)),(\overset{\bullet}{u}_k(\tau)),(u_k(x_j^{(k)})),(y_j))$,
where $\tau$ is the unit tangent vector to $D$ at $0$,
and
$$
p_{\valpha,\beta}:\CP^{\valpha,\beta}_{\vmu}\lra
\prod_k\ D^{\alpha_k}\times\CP^{\bullet\beta}_{\vmu-\valpha}
$$
sending $((u_k),((x^{(k)}_j)_k),(y_j))$ to $((u_k(0)),
(\overset{\bullet}{u}_k(\tau)),(y_j))$.

\subsubsection{} Suppose we are given $\valpha,\vbeta\in\BN[I]^K,\
\gamma\in\BN[I]$; set $\alpha:=\sum_k\ \alpha_k,\ \beta:=\sum_k\ \beta_k$.

Define a space $D^{\alpha,\beta}$ consisting of couples
$(D_{\epsilon},(x_j))$ where $D_{\epsilon}\subset D$ is some smaller disk
($0<\epsilon<1$), and $(x_j)\in D^{\alpha+\beta}$ is a configuration
such that $\alpha$ points dwell inside $D_{\epsilon}$, and $\beta$ points ---
outside $\overline{D_{\epsilon}}$.
We have an evident map
$$
q_{\alpha,\beta}:D^{\alpha,\beta}\lra D^{\alpha+\beta}.
$$

Let us define a space
$$
\CP^{\valpha,\vbeta,\gamma}_{\vmu}\subset
\tCP^K\times\prod_{k\in K}\ D^{\alpha_k,\beta_k}\times\CP^{e(\gamma)}
$$
consisting of all triples $((u_k),\bx,(y_j))$ where
$(u_k)\in\tCP^K,\ \bx\in \prod_{k}\ D^{\alpha_k,\beta_k},\
(y_j)\in\CP^{e(\gamma)}$ such that
$$
y_j\in\BP^1-\bigcup_k\ \overline{u_k(D)}.
$$
We have obvious projections
$$
q^1_{\valpha,\vbeta,\gamma}: \CP^{\valpha,\vbeta,\gamma}_{\vmu}\lra
\CP^{\valpha+\vbeta,\gamma}_{\vmu}
$$
and
$$
q^2_{\valpha,\vbeta,\gamma}: \CP^{\valpha,\vbeta,\gamma}_{\vmu}\lra
\CP^{\valpha,\beta+\gamma}_{\vmu}
$$
such that
$$
q_{\valpha+\vbeta,\gamma}\circ q^1_{\valpha,\vbeta,\gamma}=
q_{\valpha,\beta+\gamma}\circ q^2_{\valpha,\vbeta,\gamma}.
$$
We will denote the last composition by $q_{\valpha,\vbeta,\gamma}$.

We have a natural projection
$$
p_{\valpha,\vbeta,\gamma}:
\CP^{\valpha,\vbeta,\gamma}_{\vmu}\lra
\prod_k\ D_{\frac{1}{2}}^{\alpha_k}
\times\prod_k\ D_{(\frac{1}{2},1)}^{\beta_k}\times
\CP^{\bullet\gamma}_{\vmu-\valpha-\vbeta}.
$$

\subsection{}
\label{change} Let us consider a local system
$p^{\alpha*}_{\vmu}\CI^{\alpha_{\vmu}}$
over $T\CP^{o\alpha}_{\vmu}$.
By our choice of the balance function $n$, its monodromies with
respect to the rotating of tangent vectors $\omega_j$ at points $x_j$
corresponding to negative simple roots, are trivial. Therefore
it descends to a unique local system over
$\CP^{o\alpha}_{\vmu}$, to be denoted by $\CI^{\alpha}_{\mu}$.

We define a perverse sheaf
$$
\CI^{\bullet\alpha}_{\vmu}:= j_{!*}\CI^{\alpha}_{\vmu}
[\dim \CP^{\alpha}_{\vmu}]\in
\CM(\CP^{\bullet\alpha}_{\vmu};\CS).
$$

\subsection{Factorizable sheaves over $\BP^1$}
\label{fact p} Suppose we are given
a $K$-tuple of FFS's $\{\CX_k\},\ \CX_k\in\FS_{c_k},\ k\in K,\ c_k\in X/Y$,
where $\sum_k\ c_k=0$.
Let us pick
$\vmu=(\vmu_k)\geq (\lambda(\CX_k))$.

Let us call a {\em factorizable sheaf over $\BP^1$ obtained by
gluing the sheaves $\CX_k$} the following collection of data which
we will denote by $g(\{\CX_k\})$.

(i) For each $\alpha\in\BN[I]$ such that $(\vmu,\alpha)$ is admissible,
a sheaf $\CX^{\alpha}_{\vmu}\in\CM(\CP^{\alpha}_{\vmu};\CS)$.

(ii) For each $\valpha=(\alpha_k)\in\BN[I]^K,\ \beta\in \BN[I]$ such that
$(\mu,\alpha+\beta)$ is admissible (where $\alpha=\sum\ \alpha_k$),
a {\em factorization isomorphism}
$$
\phi_{\valpha,\beta}: q^*_{\valpha,\beta}\CX^{\alpha+\beta}_{\vmu}\iso
p^*_{\valpha,\beta}((\Boxtimes_{k\in K}\ \CX^{\alpha_k}_{\mu_k})
\boxtimes\CI^{\bullet\beta}_{\vmu-\valpha}).
$$

These isomorphisms should satisfy

{\em Associativity property.} The following two isomorphisms
$$
q^*_{\valpha,\vbeta,\gamma}\CX^{\alpha+\beta+\gamma}_{\vmu}\iso
p^*_{\valpha,\vbeta,\gamma}(
(\Boxtimes_k\ \CX^{\alpha_k}_{\mu_k})\boxtimes
(\Boxtimes_k\ \CI^{\bullet\beta_k}_{\mu_k-\alpha_k})\boxtimes
\CI^{\bullet\gamma}_{\vmu-\valpha-\vbeta})
$$
are equal:
$$
\psi_{\vbeta,\gamma}\circ
q^{2*}_{\valpha,\vbeta,\gamma}(\phi_{\valpha,\beta+\gamma})=
\phi_{\valpha,\vbeta}\circ
q^{1*}_{\valpha,\vbeta,\gamma}(\phi_{\valpha+\vbeta,\gamma}).
$$
Here $\psi_{\vbeta,\gamma}$ is the factorization isomorphism
for $\CI^{\bullet}$, and $\phi_{\valpha,\vbeta}$ is the tensor
product
of factorization isomorphisms for the sheaves $\CX_k$.

\subsection{Theorem}
\label{glu thm} {\em There exists a unique up to a canonical
isomorphism factorizable sheaf over $\BP^1$ obtained by gluing
the sheaves $\{\CX_k\}$.}

{\bf Proof} is similar to III.10.3. $\Box$

\newpage
\begin{center}
{\bf CHAPTER 2. Semiinfinite cohomology.}
\end{center}
\vspace{.5cm}

In this chapter we discuss, following essentially ~\cite{a},
the "Semiinfinite homological algebra" in the category $\CC$.

\section{Semiinfinite functors $Ext$ and $Tor$ in $\CC$}

\subsection{} Let us call an $\fu$-module {\em $\fu^-$-induced}
(resp., {\em $\fu^+$-induced}) if it is induced from some
$\fu^{\geq 0}$ (resp., $\fu^{\leq 0}$)-module.

\subsubsection{}
\label{tens ind} {\bf Lemma.} {\em If $M$ is a $\fu^-$-induced, and
$N$ is $\fu^+$-induced then $M\otimes_B N$ is $\fu$-projective.}

{\bf Proof.} An induced module has a filtration whose factors are
corresponding Verma modules. For Verma modules the claim is easy. $\Box$

\subsection{Definition}
\label{downup} Let $M^{\bullet}=\oplus_{\lambda\in X}\
M^{\bullet}_{\lambda}$ be a complex (possibly unbounded) in $\CC$.
We say that $M^{\bullet}$ is {\em concave} (resp. {\em convex})
if it satisfies the properties (a) and (b) below.

(a) There exists $\lambda_0\in X$ such that for any $\lambda\in X$,
if $M^{\bullet}_{\lambda}\neq 0$ then $\lambda\geq\lambda_0$
(resp. $\lambda\leq\lambda_0$).

(b) For any $\mu\in X$ the subcomplex
$\oplus_{\lambda\leq\mu}\ M^{\bullet}_{\lambda}$
(resp. $\oplus_{\lambda\leq\mu}\ M^{\bullet}_{\lambda}$)
is finite.

We will denote the category of concave (resp., convex)
complexes by $\CC^{\uparrow}$ (resp., $\CC^{\downarrow}$).

\subsection{} Let $V\in\CC$. We will say that a surjection
$\phi: P\lra V$ is {\em good} if it satisfies the following
properties:

(a) $P$ is $\fu^-$-induced;


(b) Let $\mu\in\supp P$ be an extremal point, that is, there is no
$\lambda\in\supp P$ such that $\lambda>\mu$. Then
$\mu\not\in\supp(\ker\phi)$.

For any $V$ there exists a good surjection as above. Indeed, denote by $p$
the projection $p:M(0)\lra L(0)$,
and take for $\phi$ the map $p\otimes\id_V$.

\subsection{}
\label{rel hoch}
Iteratig, we can construct a $\fu^-$-induced convex left resolution of
$B=L(0)$. Let us pick such a resolution and denote it by
$P^{\blt}_{\swarrow}$:
$$
\ldots\lra P^{-1}_{\swarrow}\lra
P^{0}_{\swarrow}\lra L(0)\lra 0
$$

We will denote by
$$
^*:\CC\lra\CC^{opp}
$$
the rigidity in $\CC$
(see e.g. ~\cite{ajs}, 7.3). We denote by $P^{\bullet}_{\nearrow}$
the complex $(P^{\bullet}_{\swarrow})^*$. It is a $\fu^-$-induced
concave right resolution of $B$. The fact that
$P^{\bullet}_{\nearrow}$ is $\fu^-$-induced follows since
$\fu^-$ is Frobenius (see e.g. ~\cite{x}).

\subsection{} In a similar manner, we can construct a
$\fu^+$-induced concave left resolution of $B$. Let us pick
such a resolution and denote it $P^{\blt}_{\nwarrow}$:
$$
\ldots\lra P^{-1}_{\nwarrow}\lra
P^{0}_{\nwarrow}\lra L(0)\lra 0
$$

We denote by $P^{\bullet}_{\searrow}$
the complex $(P^{\bullet}_{\nwarrow})^*$. It is a $\fu^+$-induced
convex right resolution of $B$.

\subsubsection{} For $M\in\CC$ we denote by
$P^{\bullet}_{\swarrow}(M)$ (resp., $P^{\bullet}_{\nearrow}(M)$,
$P^{\bullet}_{\nwarrow}(M)$, $P^{\bullet}_{\searrow}(M)$) the resolution
$P^{\bullet}_{\swarrow}\otimes_BM$ (resp.,
$P^{\bullet}_{\nearrow}\otimes_BM$, $P^{\bullet}_{\nwarrow}\otimes_BM$,
$P^{\bullet}_{\searrow}\otimes_BM$) of $M$.

\subsection{} We denote by $\CC_r$ the category of $X$-graded
right $\fu$-modules $V=\oplus_{\lambda\in X} V_{\lambda}$ such that
$$
K_i|_{V_{\lambda}}=\zeta^{-\langle i,\lambda\rangle}
$$
(note the change of a sign!), the operators $E_i$, $F_i$ acting as
$E_i:V_{\lambda}\lra V_{\lambda+i'},\ F_i:V_{\lambda}\lra V_{\lambda-i'}$.

\subsubsection{} Given $M\in\CC$, we define $M^{\vee}\in\CC_r$ as follows:
$(M^{\vee})_{\lambda}=(M_{-\lambda})^*$,
$E_i:(M^{\vee})_{\lambda}\lra (M^{\vee})_{\lambda+i'}$ is the
transpose of $E_i:M_{-\lambda-i'}\lra M_{-\lambda}$, similarly,
$F_i$ on $M^{\vee}$ is the transpose of $F_i$ on $M$.

This way we get an equivalence
$$
^{\vee}:\CC^{opp}\iso\CC_r.
$$
Similarly, one defines an equivalence $^{\vee}:\CC_r^{opp}\lra\CC$,
and we have an obvious isomorphism $^{\vee}\circ\ ^{\vee}\cong\Id$.

\subsubsection{} Given $M\in\CC$, we define $sM\in\CC_r$ as follows:
$(sM)_{\lambda}=M_{\lambda};\ xg=(sg)x$ for $x\in M, g\in\fu$ where
$$
s:\fu\lra\fu^{opp}
$$
is the antipode defined in ~\cite{ajs}, 7.2. This way we get an
equivalence
$$
s:\CC\iso\CC_r
$$
One defines an equivalence $s:\CC_r\iso\CC$ in a similar manner.
The isomorphism of functors $s\circ s\cong \Id$ is constructed
in {\em loc. cit.}, 7.3.

Note that the rigidity $*$ is just the composition
$$
^*=s\circ\ ^{\vee}.
$$

\subsubsection{} We define the categories $\CC^{\uparrow}_r$ and
$\CC^{\downarrow}_r$ in the same way as in ~\ref{downup}.

For $V\in\CC_r$ we define $P^{\bullet}_{\swarrow}(V)$ as
$$
P^{\bullet}_{\swarrow}(V)=sP^{\bullet}_{\swarrow}(sV);
$$
and $P^{\bullet}_{\nearrow}(V), P^{\bullet}_{\searrow}(V),
P^{\bullet}_{\nwarrow}(V)$ in a similar way.

\subsection{Definition} (i) Let $M,N\in\CC$. We define
$$
\Ext^{\infh+\bullet}_{\CC}(M,N):=
H^{\bullet}(\Hom_{\CC}(P^{\bullet}_{\searrow}(M),P^{\bullet}_{\nearrow}(N))).
$$

(ii) Let $V\in\CC_r,\ N\in\CC$. We define
$$
\Tor_{\infh+\bullet}^{\CC}(V,N):=
H^{-\bullet}(P^{\bullet}_{\swarrow}(V)\otimes_{\CC}P^{\bullet}_{\searrow}(N)).\
\Box
$$

Here we understand
$\Hom_{\CC}(P^{\bullet}_{\searrow}(M),P^{\bullet}_{\nearrow}(N)$ and
$P^{\bullet}_{\swarrow}(V)\otimes_{\CC}P^{\bullet}_{\searrow}(N)$ as simple
complexes associated with the corresponding double complexes.
Note that due to our boundedness properties of weights
of our resolutions, these double complexes are bounded.
Therefore all $\Ext^{\infh+i}$ and $\Tor_{\infh+i}$ spaces are finite
dimensional, and are non-zero only for finite number of $i\in\BZ$.

\subsection{Lemma.}
\label{ext tor} {\em For $M,N\in\CC$ there exist canonical
nondegenerate pairings
$$
\Ext^{\infh+n}_{\CC}(M,N)\otimes
\Tor_{\infh+n}^{\CC}(N^{\vee},M)\lra B.
$$}

{\bf Proof.} There is an evident non-degenenerate pairing
$$
\Hom_{\CC}(M,N)\otimes (N^{\vee}\otimes_{\CC}M)\lra B.
$$
It follows that the complexes computing
$\Ext$ and $\Tor$ are also canonically dual. $\Box$

\subsection{Theorem}
\label{semi thm} {\em
{\em (i)} Let $M,N\in\CC$. Let $R^{\bullet}_{\searrow}(M)$ be a
$\fu^+$-induced convex right resolution of $M$,
and $R^{\bullet}_{\nearrow}(N)$ --- a $\fu^-$-induced concave right resolution
of $N$. Then there is a canonical isomorphism
$$
\Ext^{\infh+\bullet}_{\CC}(M,N)\cong
H^{\bullet}(\Hom_{\CC}(R^{\bullet}_{\searrow}(M),R^{\bullet}_{\nearrow}(N))).
$$

{\em (ii)} Let $V\in\CC_r,N\in\CC$. Let $R^{\bullet}_{\swarrow}(V)$ be a
$\fu^-$-induced convex left resolution of $V$,
and $R^{\bullet}_{\searrow}(N)$ --- a $\fu^+$-induced convex right resolution
of $N$ lying in $\CC^{\downarrow}$.
Then there is a canonical isomorphism
$$
\Tor_{\infh+\bullet}^{\CC}(V,N)\cong
H^{-\bullet}(R^{\bullet}_{\swarrow}(V)\otimes_{\CC}R^{\bullet}_{\searrow}(N)).
$$}

{\bf Proof} will occupy the rest of the section.

\subsection{Lemma}
\label{filtr} {\em Let $V\in\CC_r$; let $R_i^{\bullet}$,
$i=1,2,$ be two $\fu^-$-induced convex left resolutions of $V$.
There exists a third $\fu^-$-induced
convex left resolution $R^{\bullet}$ of $V$, together
with two termwise surjective maps
$$
R^{\bullet}\lra R^{\bullet}_i,\ i=1,2,
$$
inducing identity on $V$.}

{\bf Proof.} We will construct $R^{\bullet}$ inductively,
from right to left. Let
$$
R^{\bullet}_i:\ \ldots\overset{d_i^{-2}}{\lra} R^{-1}_i
\overset{d_i^{-1}}{\lra} R^0_i\overset{\epsilon_i}{\lra} V.
$$
First, define $L_0:=R^0_1\times_V R^0_2$.
We denote by $\delta$ the canonical map $L^0\lra V$, and by
$q_i^0$ the projections $L^0\lra R^0_i$. Choose a good surjection
$\phi_0: R^0\lra L^0$ and define
$p^0_i:=q^0_i\circ\phi_0: R^0\lra R^0_i;\
\epsilon:=\delta\circ\phi_0:R^0\lra V$.

Set $K_i^{-1}:=\ker\epsilon_i; K:=\ker\epsilon$. The projections $p^0_i$
induce surjections $p^0_i:K^{-1}\lra K^{-1}_i$. Let us define
$$
L^{-1}:=\ker ((d_1^{-1}-p^0_1,d_2^{-1}-p_2^0):
R_1^{-1}\oplus K^{-1}\oplus R_2^{-1}\lra
K^{-1}_1\oplus K_2^{-1}).
$$
We have canonical projections
$q_i^{-1}:L^{-1}\lra R_i^{-1},\
\delta^{-1}:L^{-1}\lra K^{-1}$. Choose a good surjection
$\phi_{-1}:R^{-1}\lra L^{-1}$ and write $d^{-1}:R^{-1}\lra R^0$ for
$\delta^{-1}\circ\phi_{-1}$ composed with the inclusion
$K^{-1}\hra R^0$. We define $p_i^{-1}:=q_i^{-1}\circ\phi_{-1}$.

We have just described an induction step, and we can proceed in the same
manner. One sees directly that the left $\fu^-$-induced resolution
$R^{\bullet}$ obtained this way actually lies in $\CC^{\downarrow}_r$. $\Box$

\subsection{} Let $N\in\CC$, and let $R^{\bullet}_{\searrow}$ be
a $\fu^+$-induced convex right resolution of $N$.
For $n\geq 0$ let $b_{\geq n}(R^{\bullet}_{\searrow})$ denote the
stupid truncation:
$$
 0\lra R^0_{\searrow}\lra\ldots\lra R^n_{\searrow}\lra 0\lra\ldots
$$
For $m\geq n$ we have evident truncation maps
$b_{\geq m}(R^{\bullet}_{\searrow})\lra b_{\geq n}(R^{\bullet}_{\searrow})$.

\subsubsection{} {\bf Lemma.}
\label{stabil semi} {\em Let $R^{\bullet}_{\swarrow}$ be a
$\fu^-$-induced left resolution of a module $V\in\CC_r$. We have
$$
H^{\bullet}(R^{\bullet}_{\swarrow}\otimes_{\CC}R^{\bullet}_{\searrow})=
\invlim_n\
H^{\bullet}(R^{\bullet}_{\swarrow}\otimes_{\CC}b_{\leq n}
R^{\bullet}_{\searrow}).
$$

For every $i\in\BZ$ the inverse system
$$
\{H^i(R^{\bullet}_{\swarrow}\otimes_{\CC}b_{\leq n}R^{\bullet}_{\searrow})\}
$$
stabilizes.}

{\bf Proof.} All spaces
$H^i(R^{\bullet}_{\swarrow}\otimes_{\CC}b_{\leq n}R^{\bullet}_{\searrow})$
and only finitely many weight components of $R^{\bullet}_{\swarrow}$ and
$R^{\bullet}_{\searrow}$ contribute to $H^i$. $\Box$

\subsection{Proof of Theorem ~\ref{semi thm}} Let us consider case (ii),
and prove that $H^{\bullet}(R^{\bullet}_{\swarrow}(V)\otimes_{\CC}
R^{\bullet}_{\searrow}(N))$ does not depend, up to a canonical
isomorphism, on the choice of a resolution $R^{\bullet}_{\swarrow}(V)$.
The other independences are proved exactly in the same way.

Let $R^{\bullet}_i,\ i=1,2,$ be two left $\fu^-$-induced left convex
resolutions of $V$. According to Lemma ~\ref{filtr}, there exists a third
one, $R^{\bullet}$, projecting onto $R^{\bullet}_i$. Let us prove that
the projections induce isomorphisms
$$
H^{\bullet}(R^{\bullet}\otimes_{\CC}R^{\bullet}_{\searrow}(N))\iso
H^{\bullet}(R^{\bullet}_i\otimes_{\CC}R^{\bullet}_{\searrow}(N)).
$$
By Lemma ~\ref{stabil semi}, it suffices to prove that
$$
H^{\bullet}(R^{\bullet}\otimes_{\CC}b_{\leq n}R^{\bullet}_{\searrow}(N))\iso
H^{\bullet}(R^{\bullet}_i\otimes_{\CC}b_{\leq n}R^{\bullet}_{\searrow}(N)).
$$
for all $n$. Let $Q^{\bullet}_i$ be a cone of $R^{\bullet}\lra R^{\bullet}_i$.
It is an exact $\fu^-$-induced convex complex bounded from the right.
It is enough to check that
$H^{\bullet}(Q^{\bullet}_i\otimes_{\CC}b_{\leq n}R^{\bullet}_{\searrow}(N))=0$.

Note that for $W\in\CC_r,M\in\CC$ we have canonically
$$
W\otimes_{\CC}M=(W\otimes sM)\otimes_{\CC}B.
$$
Thus
$$
H^{\bullet}(Q^{\bullet}_i\otimes_{\CC}b_{\leq n}R^{\bullet}_{\searrow}(N))=
H^{\bullet}((Q^{\bullet}_i\otimes sb_{\leq n}R^{\bullet}_{\searrow}(N))
\otimes_{\CC}B)=0,
$$
since
$(Q^{\bullet}_i\otimes sb_{\leq n}R^{\bullet}_{\searrow}(N)$ is an exact
bounded from the right complex, consisting of modules
which are tensor products of
$\fu^+$-induced and $\fu^-$-induced, hence
$\fu$-projective modules (see Lemma ~\ref{tens ind}).

\subsubsection{} It remains to show that if $p'$ and $p''$ are two maps
between $\fu^-$-induced convex resolutions of $V$,
$R^{\bullet}_1\lra R^{\bullet}_2$,
inducing identity on $V$, then the isomorphisms
$$
H^{\bullet}(R^{\bullet}_1\otimes_{\CC}R^{\bullet}_{\searrow}(N))\iso
H^{\bullet}(R^{\bullet}_2\otimes_{\CC}R^{\bullet}_{\searrow}(N))
$$
induced by $p'$ and $p''$, coincide. Arguing as above, we see that it is enough
to prove this with $R^{\bullet}_{\searrow}(N)$ replaced by
$b_{\leq n}R^{\bullet}_{\searrow}(N)$. This in turn is equivalent to
showing that two isomorphisms
$$
H^{\bullet}((R^{\bullet}_1\otimes sb_{\leq n}R^{\bullet}_{\searrow}(N))
\otimes_{\CC}B)\iso
H^{\bullet}((R^{\bullet}_2\otimes sb_{\leq n}R^{\bullet}_{\searrow}(N))
\otimes_{\CC}B)
$$
coincide. But $R^{\bullet}_i\otimes sb_{\leq n}R^{\bullet}_{\searrow}(N)$
are complexes of projective $\fu$-modules, and the morphisms
$p'\otimes\id$ and $p''\otimes\id$ induce the same map
on cohomology, hence they are homotopic; therefore they
induce homotopic maps after tensor multiplication by $B$.

This completes the proof of the theorem. $\Box$

\section{Some calculations}

We will give a recipe for calculation of
$\Tor^{\CC}_{\infh+\blt}$ which will prove useful for the next chapter.

\subsection{} Recall that in III.13.2 the duality functor
$$
D:\CC_{\zeta^{-1}}\lra\CC^{opp}
$$
has been defined (we identify $\CC$ with $\tCC$ as usually). We will
denote objects of $\CC_{\zeta^{-1}}$ by letters with the subscript
$(\blt)_{\zeta^{-1}}$.

Note that $DL(0)_{\zeta^{-1}}=L(0)$.

Let us describe duals to Verma modules.
For $\lambda\in X$ let us denote by $M^+(\lambda)$ the
Verma module with respect to the subalgebra $\fu^+$ with the lowest
weight $\lambda$, that is
$$
M^+(\lambda):=\Ind^{\fu}_{\fu^{\leq 0}}\chi_{\lambda}
$$
where $\chi_{\lambda}$ is an evident one-dimensional representaion
of $\fu^{\leq 0}$ corresponding to the character $\lambda$.

\subsubsection{}
\label{dm} {\bf Lemma.} {\em We have
$$
DM(\lambda)_{\zeta^{-1}}=M^+(\lambda-2(l-1)\rho).
$$}

{\bf Proof} follows from ~\cite{ajs}, Lemma 4.10. $\Box$

\subsection{}
Let us denote by $K^{\blt}$ a two term complex in $\CC$
$$
L(0)\lra DM(0)_{\zeta^{-1}}
$$
concentrated in degrees $0$ and $1$, the morphism being dual to the
canonical projection $M(0)_{\zeta^{-1}}\lra L(0)_{\zeta^{-1}}$.

For $n\geq 1$ define a complex
$$
K^{\blt}_n:=b_{\geq 0}(K^{\blt\otimes n}[1]);
$$
it is concentrated in degrees from $0$ to $n-1$. For example,
$K^{\blt}_1=DM(0)_{\zeta^{-1}}$.

For $n\geq 1$ we will denote by
$$
\xi:K^{\blt}_n\lra K^{\blt}_{n+1}
$$
the map induced by the embedding $L(0)\hra DM(0)_{\zeta^{-1}}$.

We will need the following evident properties of the system
$\{K^{\blt}_n,\xi_n\}$:

(a) $K^{\blt}_n$ is $\fu^+$-induced;

(b) $K^{\blt}_n$ is exact off degrees $0$ and $n-1$;
$H^0(K^{\blt}_n)=B$. $\xi_n$ induces identity map between
$H^0(K^{\blt}_n)$ and $H^0(K^{\blt}_{n+1})$.

(c) For a fixed $\mu\in X$ there exists $m\in\BN$ such that for any
$n$ we have $(b_{\geq m}K^{\blt}_n)_{\geq \mu}=0$.
Here for $V=\oplus_{\lambda\in X}\in\CC$ we set
$$
V_{\geq\mu}:=\oplus_{\lambda\geq\mu}\ V_{\lambda}.
$$

\subsection{} Let $V\in\CC_r$; let $R^{\blt}_{\swarrow}(V)$ be a
$\fu^-$-induced convex left resolution of $V$. Let $N\in\CC$.

\subsubsection{}
\label{half tor}{\bf Lemma.} (i) {\em For a fixed $k\in\BZ$ the direct system
$\{H^k(R^{\blt}_{\swarrow}(V)\otimes_{\CC}(K^{\blt}_n\otimes N)),\xi_n\}$
stabilizes.}

(ii) {\em We have a canonical isomorphism
$$
\Tor^{\CC}_{\infh+\blt}(V,N)\cong
\dirlim_n\ H^{-\blt}(R^{\blt}_{\swarrow}(V)\otimes_{\CC}(K^{\blt}_n\otimes N)).
$$}

{\bf Proof.} (i) is similar to Lemma ~\ref{stabil semi}.
(ii) By Theorem ~\ref{semi thm} we can use any $\fu^+$-induced
right convex resolution of $N$ to compute
$\Tor^{\CC}_{\infh+\blt}(V,N)$. Now extend
$K^{\blt}_n\otimes N$ to a $\fu^+$-induced convex resolution of $N$
and argue like in the proof of Lemma ~\ref{stabil semi} again. $\Box$

\subsection{} Recall the notations of ~\ref{rel hoch}
and take for $R^{\blt}_{\swarrow}(V)$ the resolution
$P^{\blt}_{\swarrow}(V)=P^{\blt}_{\swarrow}\otimes V$. Then
$$
H^{\blt}(P^{\blt}_{\swarrow}(V)\otimes_{\CC}(K^{\blt}_n\otimes N))=
H^{\blt}(V\otimes_{\CC}(P^{\blt}_{\swarrow}\otimes K^{\blt}_n
\otimes N)).
$$
Note that $P^{\blt}_{\swarrow}\otimes K^{\blt}_n\otimes N$ is a right
bounded complex quasi-isomorphic to $K^{\blt}_n\otimes N$.
The terms of $P^{\blt}_{\swarrow}\otimes K^{\blt}_n$ are
$\fu$-projective by Lemma ~\ref{tens ind}, hence the terms of
$P^{\blt}_{\swarrow}\otimes K^{\blt}_n\otimes N$ are projective by rigidity
of $\CC$. Therefore,
$$
H^{-\blt}(V\otimes_{\CC}(P^{\blt}_{\swarrow}\otimes K^{\blt}_n
\otimes N))=
\Tor^{\CC}_{\blt}(V,K^{\blt}_n\otimes N).
$$
Here $\Tor^{\CC}_{\blt}(*,*)$ stands for the zeroth
weight component of  $\Tor^{\fu}_{\blt}(*,*)$.

Putting all the above together, we get

\subsection{Corollary}
{\em For a fixed $k\in\BZ$ the direct system
$\{\Tor_k^{\CC}(V,K^{\blt}_n\otimes N)\}$ stabilizes.

We have
$$
\Tor^{\CC}_{\infh+\blt}(V,N)=
\dirlim_n\ \Tor^{\CC}_{\blt}(V,K^{\blt}_n\otimes N).\ \Box
$$}

\subsection{} Dually, consider complexes $DK^{\blt}_{n,\zeta^{-1}}$. They form
a projective system
$$
\{\ldots\lra DK^{\blt}_{n+1,\zeta^{-1}}\lra DK^{\blt}_{n,\zeta^{-1}}\lra
\ldots \}
$$

These complexes enjoy properties dual to (a) --- (c) above.

\subsection{Theorem}
\label{two side} {\em For every $k\in\BZ$ we have canonical
isomorphisms
$$
\Tor^{\CC}_{\infh+k}(V,N)\cong
\invlim_m\dirlim_n\ H^{-k}((V\otimes sDK^{\blt}_{m,\zeta^{-1}})\otimes_{\CC}
(K^{\blt}_n\otimes N)).
$$
Both the inverse and the direct systems actually stabilize.}

{\bf Proof} follows from Lemma ~\ref{half tor}. We leave details
to the reader. $\Box$

Here is an example of calculation of $\Tor^{\CC}_{\infh+\blt}$.

\subsection{Lemma}
\label{tor rho} {\em
$\Tor^{\CC}_{\infh+\blt}(B,L(2(l-1)\rho))=B\ \mbox{in degree $0$.}$}

{\bf Proof.} According to Lemma ~\ref{ext tor} it suffices to prove
that
$\Ext_{\CC}^{\infh+\blt}(L(2(l-1)\rho),L(0))=B$. Choose a $\fu^+$-induced
right convex resolution
$$
L(2(l-1)\rho)\overset{\epsilon}{\lra} R^{\blt}_{\searrow}
$$
such that
$$
R^0_{\searrow}= DM(2(l-1)\rho)_{\zeta^{-1}}=M^+(0),
$$
and all the weights in $R^{\geq 1}_{\searrow}$ are $<2(l-1)\rho$.

Similarly, choose a $\fu^-$-free right concave resolution
$$
L(0)\overset{\epsilon}{\lra} R^{\blt}_{\nearrow}
$$
such that
$$
R^0_{\nearrow}=M(2(l-1)\rho)=DM^+(0)_{\zeta^{-1}},
$$
and all the weights of $R^{\geq 1}_{\nearrow}$ are $>0$. By
Theorem ~\ref{semi thm} we have
$$
\Ext_{\CC}^{\infh+\blt}(L(2(l-1)\rho),L(0))=
H^{\blt}(\Hom_{\CC}^{\blt}(R^{\blt}_{\searrow},R^{\blt}_{\nearrow})).
$$
Therefore it is enough to prove that

(a) $\Hom_{\CC}(R^0_{\searrow},R^0_{\nearrow})=B$;

(b) $\Hom_{\CC}(R^m_{\searrow},R^n_{\nearrow})=0$ for $(m,n)\neq (0,0)$.

(a) is evident. Let us prove (b) for, say, $n>0$. $R^m_{\searrow}$ has
a filtration with successive quotients of type
$M^+(\lambda)$, $\lambda\leq 0$; similarly,
$R^n_{\searrow}$ has
a filtration with successive quotients of type
$DM^+(\mu)_{\zeta^{-1}}$, $\mu > 0$. We have
$\Hom_{\CC}(M^+(\lambda),DM^+(\mu)_{\zeta^{-1}})=0$, therefore
$\Hom_{\CC}(R^m_{\searrow},R^n_{\nearrow})=0$. The proof for $m>0$
is similar. Lemma is proven. $\Box$

\vspace{.5cm}
{\em CONFORMAL BLOCKS AND $Tor^{\CC}_{\infh+\blt}$}
\vspace{.5cm}

\subsection{}
\label{triv sub} Let $M\in \CC$. We have a canonical embedding
$$
\Hom_{\CC}(B,M)\hra M
$$
which identifies $\Hom_{\CC}(B,M)$ with the maximal trivial
subobject of $M$. Dually, we have a canonical
epimorhism
$$
M\lra\Hom_{\CC}(M,B)^*
$$
which identifies $\Hom_{\CC}(M,B)^*$ with the maximal trivial quotient
of $M$. Let us denote by $\langle M\rangle$ the image
of the composition
$$
\Hom_{\CC}(B,M)\lra M\lra\Hom_{\CC}(M,B)^*
$$
Thus, $\langle M\rangle$ is canonically a subquotient of $M$.

One sees easily that if $N\subset M$ is a trivial direct summand of $M$
which is maximal, i.e. not contained in greater direct summand, then
we have a canonical isomorphism $\langle M\rangle\iso N$. By this reason,
we will call $\langle M\rangle$ {\em the maximal trivial
direct summand} of $M$.

\subsection{}
\label{alc} Let
$$
\Delta=\{\lambda\in X|\ \langle i,\lambda+\rho\rangle>0,\
\mbox{for all }i\in I;\ \langle\gamma,\lambda+\rho\rangle<l\}
$$
denote the first alcove. Here $\gamma\in \CR\subset Y$ is the
highest coroot.

For $\lambda_1,\ldots,\lambda_n\in\Delta$ {\em the space of conformal
blocks} is defined as
$$
\langle L(\lambda_1),\ldots,L(\lambda_n)\rangle :=
\langle L(\lambda_1)\otimes \ldots\otimes L(\lambda_n)\rangle
$$
(see e.g. ~\cite{an} and Lemma ~\ref{compar conf} below).

\subsection{Corollary}
\label{conf subq} {\em The space of conformal blocks
$\langle L(\lambda_1),\ldots,L(\lambda_n)\rangle$
is canonically a subquotient of
$\Tor^{\CC}_{\infh+0}(B,L(\lambda_1)\otimes \ldots\otimes
L(\lambda_n)\otimes L(2(l-1)\rho)).$}

{\bf Proof} follows easily from the definition of $\langle\blt\rangle$
and Lemma ~\ref{tor rho}. $\Box$

\subsection{}
\label{ex conf} Let us consider an example showing that
$\langle L(\lambda_1),\ldots,L(\lambda_n)\rangle$ is in general a
{\em proper} subquotient of
$\Tor^{\CC}_{\infh+0}(B,L(\lambda_1)\otimes\ldots\otimes
L(\lambda_n)\otimes L(2(l-1)\rho))$.

We leave the following to the reader.

\subsubsection{} {\bf Exercise.} {\em Let $P(0)$ be the indecomposable
projective cover of $L(0)$. We have
$\Tor^{\CC}_{\infh+0}(B,P(0))=B.\ \Box$}

We will construct an example featuring $P(0)$ as a direct summand of
$L(\lambda_1)\otimes \ldots\otimes L(\lambda_n)$.

Let us take a root datum of type $sl(2)$; take $l=5,\ n=4,\
\lambda_1=\lambda_2=2,\ \lambda_3=\lambda_4=3$
(we have identified $X$ with $\BZ$).

In our case $\rho=1$, so $2(l-1)\rho=8$. Note that $P(0)$ has highest
weight $8$, and it is a unique indecomposable projective with
the highest weight $8$. So, if we are able to find a
projective summand of highest weight $0$ in $V=L(2)\otimes L(2)\otimes
L(3)\otimes L(3)$ then $V\otimes L(8)$ will
contain a projective summand of highest weight $8$, i.e. $P(0)$.

Let $U_B$ denote the quantum group with divided powers over $B$
(see ~\cite{l2}, 8.1). The algebra $\fu$ lies inside $U_B$. It is wellknown
that all irreducibles $L(\lambda)$, $\lambda\in\Delta$, lift to simple
$U_B$-modules $\wt{L(\lambda)}$ and for
$\lambda_1,\ldots,\lambda_n\in\Delta$ the $U_B$-module
$\wt{L(\lambda_1)}\otimes\ldots\otimes\wt{L(\lambda_n)}$ is a direct sum
of irreducibles $\wt{L(\lambda)},\ \lambda\in\Delta,$
and indecomposable projectives $\wt{P(\lambda)},\ \lambda\geq 0$
(see, e.g. ~\cite{an}).

Thus $\wt{L(2)}\otimes\wt{L(2)}\otimes\wt{L(3)}\otimes\wt{L(3)}$
contains an indecomposable projective summand with the highest weight
$10$, i.e. $\wt{P(8)}$. One can check easily
that when restricted to $\fu$, $\wt{P(8)}$ remains projective
and contains a summand $P(-2)$. But the highest weight of $P(-2)$ is zero.

We conclude that
$L(2)\otimes L(2)\otimes L(3)\otimes L(3)\otimes L(8)$ contains
a projective summand $P(0)$, whence
$$
\langle L(2),L(2),L(3),L(3)\rangle\neq
\Tor^{\CC}_{\infh+0}(B,L(2)\otimes L(2)\otimes L(3)\otimes L(3)\otimes L(8)).
$$

\newpage
\begin{center}
{\bf CHAPTER 3. Global sections.}
\end{center}
\vspace{.5cm}

\section{Braiding and balance in $\CC$ and $\FS$}

\subsection{}
\label{br} Let $U_B$ be the quantum group with divided powers,
cf ~\cite{l2}, 8.1. Let $_R\CC$ be the category of finite dimensional
integrable $U_B$-modules defined in ~\cite{kl}IV, \S 37.
It is a rigid braided tensor category. The braiding, i.e.
family of isomorphisms
$$
\tR_{V,W}:V\otimes W\iso W\otimes V,\ V,W\in\ _R\CC,
$$
satisfying the usual constraints, has been defined in
{}~\cite{l}, Ch. 32.


\subsection{}
\label{ups} As $\fu$ is a subalgebra of $U_B$, we have the restriction
functor preserving $X$-grading
$$
\Upsilon:\ _R\CC\lra\CC.
$$
The following theorem is due to G.Lusztig (private communication).

\subsubsection{}
\label{braid thm} {\bf Theorem.} {\em
{\em (a)} There is a unique braided structure $(R_{V,W},\theta_V)$
on $\CC$ such that the restriction functor $\Upsilon$ commutes
with braiding.

{\em (b)} Let $V=L(\lambda)$, and let $\mu$ be the highest weight of
$W\in\CC$, i.e. $W_{\mu}\neq 0$ and $W_{\nu}\neq 0$ implies
$\nu\leq\mu$. Let $x\in V,\ y\in W_{\mu}$. Then
$$
R_{V,W}(x\otimes y)=\zeta^{\lambda\cdot\mu}y\otimes x;
$$

{\em (c)} Any braided structure on $\CC$
enjoying the property {\em (b)} above coincides with that defined in
{\em (a)}.
\ $\Box$}

\subsection{} Recall that an automorphism $\ttheta=\{\ttheta_V: V\iso V\}$
of the identity functor of $_R\CC$ is called {\em balance}
if for any $V,W\in\ _R\CC$ we have
$$
\tR_{W,V}\circ\tR_{V,W}=\ttheta_{V\otimes W}\circ (\ttheta_V\otimes
\ttheta_W)^{-1}.
$$

The following proposition is an easy application of the results
of ~\cite{l}, Chapter 32.

\subsection{Proposition} {\em The category $_R\CC$ admits a unique balance
$\ttheta$ such that

--- if $\tL(\lambda)$ is an irreducible in $_R\CC$ with the
highest weight $\lambda$, then $\ttheta$ acts on $\tL(\lambda)$ as
multiplication by $\zeta^{n(\lambda)}$.\ $\Box$}

Here $n(\lambda)$ denotes the function introduced in ~\ref{bal fun}.

Similarly to ~\ref{braid thm}, one can prove

\subsection{Theorem} {\em {\em (a)} There is a unique balance
$\theta$ on $\CC$ such that $\Upsilon$ commutes with balance;

{\em (b)} $\theta_{L(\lambda)}=\zeta^{n(\lambda)}$;

{\em (c)} If $\theta'$ is a balance in $\CC$ having property {\em (b)},
then $\theta'=\theta$.\ $\Box$}

\subsection{}
\label{phi braid} According to Deligne's ideology, ~\cite{d1}, the gluing
construction of ~\ref{glu thm} provides the category $\FS$ with the
balance $\theta^{\FS}$. Recall that the braiding $R^{\FS}$ has
been defined in III.11.11 (see also 11.4). It follows easily from the
definitions that $(\Phi(R^{\FS}),\Phi(\theta^{\FS}))$ satisfy the
properties (b)(i) and (ii) above. Therefore, we have
$(\Phi(R^{\FS}),\Phi(\theta^{\FS}))=(R,\theta)$, i.e.
$\Phi$ is an equivalence of braided balanced categories.

\section{Global sections over $\CA(K)$}
\label{global}

\subsection{} Let $K$ be a finite non-empty set, $|K|=n$, and let
$\{\CX_k\}$ be a $K$-tuple of finite gactorizable sheaves. Let $\lambda_k:=
\lambda(\CX_k)$ and $\lambda=\sum_k\ \lambda_k$; let $\alpha\in\BN[I]$.
Consider the sheaf $\CX^{\alpha}(K)$ over $\CA^{\alpha}_{\vlambda}(K)$
obtained by gluing $\{\CX_k\}$, cf. III.10.3. Thus
$$
\CX^{\alpha}(K)=g_K(\{\CX_k\})^{\alpha}_{\vlambda}
$$
in the notations of {\em loc.cit.}

We will denote by $\eta$, or sometimes by $\eta^{\alpha}$,
or $\eta^{\alpha}_K$ the projection
$\CA^{\alpha}(K)\lra\CO(K)$. We are going to describe
$R\eta_*\CX^{\alpha}(K)[-n]$. Note that it is an element of
$\CD(\CO(K))$  which is smooth, i.e. its cohomology sheaves are local
systems.

\subsection{} Let $V_1,\ldots,V_n\in\CC$. Recall (see II.3) that
$C^{\blt}_{\fu^-}(V_1\otimes\ldots\otimes V_n)$ denotes the
Hochschild complex of the $\fu^-$-module $V_1\otimes\ldots\otimes V_n$. It
is naturally $X$-graded, and its $\lambda$-component is denoted
by the subscript $(\blt)_{\lambda}$ as usually.

Let us consider a homotopy point $\bz=(z_1,\ldots,z_n)\in\CO(K)$
where all $z_i$ are real, $z_1<\ldots<z_n$. Choose a bijection
$K\iso [n]$. We want to describe a stalk
$R\eta_*\CX^{\alpha}(K)_{\bz}[-n]$. The following theorem generalizes
Theorem II.8.23. The proof is similar to {\em loc. cit}, cf. III.12.16,
and will appear later.

\subsection{Theorem}
\label{coh tens} {\em There is a canonical isomorphism,
natural in $\CX_i$,
$$
R\eta_*\CX^{\alpha}(K)_{\bz}[-n]\cong
C^{\blt}_{\fu^-}(\Phi(\CX_1)\otimes\ldots\otimes\Phi(\CX_n))_{\lambda-\alpha}.
\ \Box
$$}

\subsection{} The group $\pi_1(\CO(K);\bz)$ is generated
by counterclockwise loops of $z_{k+1}$ around $z_k$, $\sigma_i$,
$k=1,\ldots,n-1$. Let $\sigma_k$ act on
$\Phi(\CX_1)\otimes\ldots\otimes\Phi(\CX_n)$ as
$$
\id\otimes\ldots\otimes R_{\Phi(\CX_{k+1}),\Phi(\CX_k)}\circ
R_{\Phi(\CX_{k}),\Phi(\CX_{k+1})}\otimes\ldots\otimes\id.
$$
This defines an action of $\pi_1(\CO(K);\bz)$ on
$\Phi(\CX_1)\otimes\ldots\otimes\Phi(\CX_n)$,
whence we get an action of this group on
$C^{\blt}_{\fu^-}(\Phi(\CX_1)\otimes\ldots\otimes\Phi(\CX_n))$ respecting
the $X$-grading. Therefore we get a complex of local systems
over $\CO(K)$; let us denote it
$C^{\blt}_{\fu^-}(\Phi(\CX_1)\otimes\ldots\otimes\Phi(\CX_n))^{\heartsuit}$.

\subsection{Theorem} {\em There is a canonical isomorphism
in $\CD(\CO(K))$
$$
R\eta_*\CX^{\alpha}(K)[-n]\iso
C^{\blt}_{\fu^-}(\Phi(\CX_1)\otimes\ldots\otimes\Phi(\CX_n))^{\heartsuit}_{\lambda
-\alpha}.
$$}

{\bf Proof} follows from ~\ref{phi braid} and  Theorem ~\ref{coh tens}.
$\Box$

\subsection{Corollary} {\em Set
$\lambda_{\infty}:=\alpha+2(l-1)\rho-\lambda.$
There is a canonical isomorphism in $\CD(\CO(K))$
$$
R\eta_*\CX^{\alpha}(K)[-n]\iso
C^{\blt}_{\fu}(\Phi(\CX_1)\otimes\ldots\otimes\Phi(\CX_n)
\otimes DM(\lambda_{\infty})_{\zeta^{-1}})^{\heartsuit}_{0}.
$$}

{\bf Proof.} By Shapiro's lemma, we have a canonical morphism
of complexes which is a quasiisomorphism
$$
C^{\blt}_{\fu^-}(\Phi(\CX_1)\otimes\ldots\otimes\Phi(\CX_n))_{\lambda
-\alpha}\lra
C^{\blt}_{\fu}(\Phi(\CX_1)\otimes\ldots\otimes\Phi(\CX_n)
\otimes M^+(\alpha-\lambda))_{0}.
$$
By Lemma ~\ref{dm}, $M^+(\alpha-\lambda)=
DM(\lambda_{\infty})_{\zeta^{-1}}$. $\Box$

\section{Global sections over $\CP$}

\subsection{} Let $J$ be a finite set, $|J|=m$, and
$\{\CX_j\}$ a $J$-tuple of finite factorizable sheaves. Set
$\mu_j:=\lambda(\CX_j),\ \vmu=(\mu_j)\in X^J$. Let
$\alpha\in \BN[I]$ be such that $(\vmu,\alpha)$ is admissible,
cf. ~\ref{bal fun}. Let $\CX^{\alpha}_{\vmu}$ be the preverse
sheaf on $\CP^{\alpha}_{\vmu}$ obtained by gluing the sheaves $\CX_j$,
cf. ~\ref{fact p} and ~\ref{glu thm}.

Note that the group $\PGL_2(\BC)=\Aut(\BP^1)$ operates naturally on
$\CP^{\alpha}_{\vmu}$ and the sheaf $\CX^{\alpha}_{\vmu}$ is equivariant
with respect to this action.

Let
$$
\breta:\CP^{\alpha}_{\vmu}\lra T\CP^{oJ}
$$
denote the natural projection; we will denote this map also by
$\breta_J$ or $\breta_J^{\alpha}$. Note that $\breta$ commutes
with the natural action of $\PGL_2(\BC)$ on these spaces.
Therefore $R\breta_*\CX^{\alpha}_{\vmu}$ is a smooth
$\PGL_2(\BC)$-equivariant
complex on $T\CP^{oJ}$. Our aim in this section will be to compute
this complex algebraically.

Note that $R\breta_*\CX^{\alpha}_{\vmu}$ descends uniquely
to the quotient
$$
\ul{T\CP}^{oJ}:= T\CP^{oJ}/\PGL_2(\BC)
$$

\subsection{} Let us pick a bijection $J\iso [m]$.
Let $\ul{\CZ}$ be a contractible real submanifold
of $\ul{T\CP}^{oJ}$ defined in ~\cite{kl}II, 13.1
(under the name $\ul{\CV}_0$). Its points are configurations
$(z_1,\tau_1,\ldots,z_m,\tau_m)$ such that $z_j\in\BP^1(\BR)=S\subset
\BP^1(\BC)$; the points $z_j$ lie on $S$ in this cyclic order;
they orient $S$ in the same way as $(0,1,\infty)$ does;
the tangent vectors $\tau_j$ are real and compatible with this orientation.

\subsection{Definition}
\label{pos} An $m$-tuple of weights $\vmu\in X^m$
is called {\em positive} if
$$
\sum_{j=1}^m\ \mu_j+(1-l)2\rho\in\BN[I]\subset X.
$$
If this is so, we will denote
$$
\alpha(\vmu):= \sum_{j=1}^m\ \mu_j+(1-l)2\rho\ \Box
$$

\subsection{Theorem}
\label{stalk glob} {\em Let $\CX_1,\ldots,\CX_m\in\FS$.
Let $\vmu$ be a positive $m$-tuple of weights, $\mu_j\geq\lambda(\CX_j)$,
and let $\alpha=\alpha(\vmu)$. Let $\CX^{\alpha}_{\vmu}$ be the sheaf
on $\CP^{\alpha}_{\vmu}$ obtained by gluing the sheaves $\CX_j$. There is a
canonical isomorphism
$$
R^{\blt}\breta_*\CX^{\alpha}_{\vmu}[-2m]_{\CZ}\cong
\Tor^{\CC}_{\infh-\blt}(B,\Phi(\CX_1)\otimes\ldots\otimes\Phi(\CX_m)).
$$}

{\bf Proof} is sketched in the next few subsections.

\subsection{Two-sided \v{C}ech resolutions} The idea of the construction
below is inspired by ~\cite{b}, p. 40.

Let $P$ be a topological space, $\CU=\{U_i|\ i=1,\ldots, N\}$ an open
covering of $P$. Let $j_{i_0i_1\ldots i_a}$ denote the embedding
$$
U_{i_0}\cap\ldots\cap U_{i_a}\hra P.
$$

Given a sheaf $\CF$ on $P$, we have a canonical morphism
\begin{equation}
\label{right ch}
\CF\lra\Cch^{\blt}(\CU;\CF)
\end{equation}
where
$$
\Cch^a(\CU;\CF)=\oplus_{i_0<i_1<\ldots<i_a}\ j_{i_0i_1\ldots i_a*}
j_{i_0i_1\ldots i_a}^*\CF,
$$
the differential being the usual \v{C}ech one.

Dually, we define a morphism
\begin{equation}
\label{left ch}
\Cch_{\blt}(\CU;\CF)\lra\CF
\end{equation}
where
$$
\Cch_{\blt}(\CU;\CF):\
0\lra \Cch_{N}(\CU;\CF)\lra \Cch_{N-1}(\CU;\CF)\lra\ldots\lra
\Cch_{0}(\CU;\CF)\lra 0
$$
where
$$
\Cch_a(\CU;\CF)=\oplus_{i_0<i_1<\ldots<i_a}\ j_{i_0i_1\ldots i_a!}
j_{i_0i_1\ldots i_a}^*\CF.
$$
If $\CF$ is injective then the arrows ~(\ref{right ch}) and
{}~(\ref{left ch}) are quasiisomorphisms.

Suppose we have a second open covering of $P$, $\CV=\{V_j|\
j=1,\ldots, N\}$. Let us define sheaves
$$
\Cch^a_b(\CU,\CV;\CF):=
\Cch_b(\CV;\Cch^a(\CU;\CF));
$$
they form a bicomplex. Let us consider the associated simple complex
$\Cch^{\blt}(\CU,\CV;\CF)$, i.e.
$$
\Cch^{i}(\CU,\CV;\CF)=\oplus_{a-b=i}\ \Cch^a_b(\CU,\CV;\CF).
$$
It is a complex concentrated in degrees from $-N$ to $N$. We have
canonical morphisms
$$
\CF\lra\Cch^{\blt}(\CU;\CF)\lla\Cch^{\blt}(\CU,\CV;\CF)
$$
If $\CF$ is injective then both arrows are quasiisomorphisms, and
the above functors are exact on injective sheaves. Therefore,
they pass to derived categories, and we get a functor
$\CK\mapsto\Cch^{\blt}(\CU,\CV;\CK)$ from the bounded derived
category $\CD^b(P)$ to the bounded filtered derived category
$\CD F(P)$. This implies

\subsection{Lemma}
\label{bi cech} {\em Suppose that $\CK\in\CD^b(P)$ is
such that $R^i\Gamma(P;\Cch^a_b(\CU,\CV;\CK))
=0$ for all $a,b$ and all $i\neq 0$. Then we have a canonical isomorphism
in $\CD^b(P)$,
$$
R\Gamma(P;\CK)\iso R^0\Gamma(P;\Cch^{\blt}(\CU,\CV;\CK)).\ \Box
$$}

\subsection{} Returning to the assumptions of
theorem ~\ref{stalk glob}, let us pick a point
$\bz=(z_1,\tau_1,\ldots,z_m,\tau_m)\in T\CP^{oJ}$ such that $z_j$ are
real numbers
$z_1<\ldots <z_m$ and tangent vectors are directed to the right.

By definition, we have canonically
$$
R\breta_*(\CP^{\alpha}_{\vmu};\CX^{\alpha}_{\vmu})_{\ul{\CZ}}=
R\Gamma(\CP^{\alpha};\CK)
$$
where
$$
\CK:=\CX^{\alpha}_{\vmu}|_{\breta^{-1}(\bz)}[-2m].
$$
Let us pick $N\geq |\alpha|$ and reals $p_1,\ldots, p_N,q_1,\ldots,q_N$
such that
$$
p_1<\ldots <p_N<z_1<\ldots <z_m<q_N<\ldots<q_1.
$$
Let us define two open coverings $\CU=\{U_i|\ i=1,\ldots,N\}$ and
$\CV=\{V_i|\ i=1,\ldots, N\}$ of the space $\CP^{\alpha}$ where
$$
U_i=\CP^{\alpha}-\bigcup_k\ \{t_k=p_i\};\
V_i=\CP^{\alpha}-\bigcup_k\ \{t_k=q_i\},
$$
where $t_k$ denote the standard coordinates.

\subsection{Lemma} (i) {\em We have
$$
R^i\Gamma(\CP^{\alpha};\Cch^a_b(\CU,\CV;\CK))=0
$$
for all $a,b$ and all $i\neq 0$.}

(ii) {\em We have canonical isomorphism
$$
R^0\Gamma(\CP^{\alpha};\Cch^{\blt}(\CU,\CV;\CK))\cong
sDK^{\blt}_{N,\zeta^{-1}}\otimes_{\CC}
(K^{\blt}_N\otimes\Phi(\CX_1)\otimes\ldots\otimes\Phi(\CX_m)),
$$
in the notations of ~\ref{two side}.}

{\bf Proof} (sketch).  We should regard the computation
of $R\Gamma(\CP^{\alpha};\Cch^a_b(\CU,\CV;\CK))$ as the computation
of global sections over $\CP^{\alpha}$ of a sheaf obtained
by gluing $\CX_j$ into points $z_j$,
the Verma sheaves $\CM(0)$ or irreducibles $\CL(0)$ into the points
$p_j$, and dual sheaves $D\CM(0)_{\zeta^{-1}}$ or
$D\CL(0)_{\zeta^{-1}}$ into the points $q_j$.

Using $\PGL_2(\BR)$-invariance, we can
move one of the points $p_j$ to infinity. Then, the desired global sections
are reduced to global sections over an affine space $\CA^{\alpha}$, which are
calculated by means of Theorem ~\ref{coh tens}.

Note that in our situation all the sheaves $\Cch^a_b(\CU,\CV;\CK)$ actually
belong to the abelian category $\CM(\CP^{\alpha})$ of perverse
sheaves. So $\Cch^{\blt}(\CU,\CV;\CK)$ is a resolution of $\CK$ in
$\CM(\CP^{\alpha})$.
$\Box$

\subsection{} The conclusion of ~\ref{stalk glob} follow from
the previous lemma and Theorem ~\ref{two side}.
$\Box$

\subsection{} The group $\pi_1(\ul{T\CP}^{om},\CZ)$ operates on
the spaces
$\Tor^{\CC}_{\infh+\blt}(B,\Phi(\CX_1)\otimes\ldots\otimes\Phi(\CX_m))$
via its action on the object $\Phi(\CX_1)\otimes\ldots\otimes\Phi(\CX_m)$
induced by the braiding and balance in $\CC$. Let us denote by
$$
 \Tor^{\CC}_{\infh+\blt}(B,\Phi(\CX_1)\otimes\ldots\otimes\Phi(\CX_m))
^{\heartsuit}
$$
the corresponding local system on $\ul{T\CP}^{om}$.

\subsection{Theorem}
\label{global thm} {\em There is a canonical isomorphism of local
systems on $\ul{T\CP}^{om}$:
$$
R^{\blt-2m}\breta_*\CX^{\alpha}_{\vmu}\cong
\Tor^{\CC}_{\infh-\blt}(B,\Phi(\CX_1)\otimes\ldots\otimes\Phi(\CX_m))
^{\heartsuit}
$$}

{\bf Proof} follows immediately from ~\ref{phi braid} and
Theorem ~\ref{stalk glob}. $\Box$

\section{Application to conformal blocks}

\subsection{} In applications to conformal blocks we will encounter
the roots of unity $\zeta$ of not necessarily odd degree $l$. So we have
to generalize all the above considerations to the case of arbitrary $l$.

The definitions of the categories $\CC$ and $\FS$ do not change (for
the category $\CC$ the reader may consult ~\cite{ap}, \S3). The construction
of the functor $\Phi: \FS\lra\CC$ and the proof that $\Phi$ is an
equivalence repeats the one in III word for word.

Here we list the only minor changes (say, in the definition of the Steinberg
module) following ~\cite{l} and ~\cite{ap}.

\subsubsection{}
So suppose $\zeta$ is a primitive root of unity of an {\em even} degree $l$.

We define $\ell:=\frac{l}{2}$. For the sake of unification of notations,
in case $l$ is {\em odd} we define $\ell:=l$.
For $i\in I$ we define $\ell_i:=\frac{\ell}
{(\ell,d_i)}$ where $(\ell,d_i)$ stands for the greatest common divisor
of $\ell$ and $d_i$.

For a coroot $\alpha\in\CR\in Y$ we can find an element $w$ of the Weyl
group $W$ and a simple coroot $i\in Y$ such that $w(i)=\alpha$ (notations
of ~\cite{l}, 2.3). We define $\ell_\alpha:=\frac{\ell}{(\ell,d_i)}$, and
the result does not depend on a choice of $i$ and $w$.

We define $\gamma_0\in\CR$ to be the highest coroot, and $\beta_0\in\CR$
to be the coroot dual to the highest root. Note that $\gamma_0=\beta_0$ iff
our root datum is simply laced.

\subsubsection{}
We define
$$
Y^*_{\ell}:=\{\lambda\in X|\lambda\cdot\mu\in\ \ell\BZ\ \mbox{for any }\mu\in
X\}
$$

One should replace the congruence modulo $lY$ in the Definition ~\ref{admis}
and in ~\ref{bal fun} by the congruence modulo $Y^*_{\ell}$.

We define $\rho_\ell\in X$ as the unique element such that
$\langle i,\rho_\ell\rangle=\ell_i-1$ for any $i\in I$.

Then the {\em Steinberg module} $L(\rho_\ell)$ is irreducible projective
in $\CC$ (see ~\cite{ap} 3.14).

Note also that $\rho_\ell$ is the highest weight of $\fu^+$.

One has to replace
all the occurences of $(l-1)2\rho$ in the above sections
by $2\rho_\ell$.

In particular, the new formulations of the Definition ~\ref{pos} and
the Theorem ~\ref{stalk glob} force us to make the following changes
in ~\ref{bal fun} and ~\ref{change}.

In ~\ref{bal fun} we choose a balance function $n$ in the form
$$n(\mu)=\frac{1}{2}\mu\cdot\mu-\mu\cdot\rho_\ell$$
In other words, we set $\nu_0=-\rho_\ell$.
{\em This balance function does not necessarily have the property that
$n(-i')\equiv 0 \bmod l$. It is only true that $n(-i')\equiv 0\bmod\ell$.}

We say that a pair $(\vmu,\alpha)$ is admissible if $\sum_k\mu_k-\alpha
\equiv 2\rho_\ell\bmod Y^*_{\ell}$.


\subsubsection{}
\label{alcove}
The last change concerns the definition of the first alcove in ~\ref{alc}.

The corrected definition reads as follows:

if $\ell_i=\ell$ for any $i\in I$, then
$$
\Delta_l=\{\lambda\in X|\ \langle i,\lambda+\rho\rangle>0,\
\mbox{for all }i\in I;\ \langle\gamma_0,\lambda+\rho\rangle<\ell\};
$$
if not, then
$$
\Delta_l=\{\lambda\in X|\ \langle i,\lambda+\rho\rangle>0,\
\mbox{for all }i\in I;\ \langle\beta_0,\lambda+\rho\rangle<\ell_{\beta_0}\}
$$

\subsection{}
Let $\hfg$ denote the affine
Lie algebra associated with $\fg$:
$$
0\lra\BC\lra\hfg\lra\fg((\epsilon))\lra 0.
$$
Let $\tCO_\kappa$ be the category of integrable $\hfg$-modules with the central
charge $\kappa-\check{h}$ where $\check{h}$ stands for the dual
 Coxeter number of $\fg$.
It is a semisimple balanced braided rigid tensor category
(see e.g. ~\cite{ms} or ~\cite{f2}).

Let $\CO_{-\kappa}$ be the category of $\fg$-integrable $\hfg$-modules
of finite length with the central charge $-\kappa-\check{h}$. It is a
balanced braided rigid tensor (bbrt) category (see ~\cite{kl}).
Let $\tCO_{-\kappa}$ be the semisimple subcategory of $\CO_{-\kappa}$
formed by direct sums of simple $\hfg$-modules with highest weights
in the alcove $\nabla_\kappa$:
$$
\nabla_\kappa:=\{\lambda\in X|\ \langle i,\lambda+\rho\rangle>0,\
\mbox{for all }i\in I;\ \langle\beta_0,\lambda+\rho\rangle<\kappa\}
$$
The bbrt structure
on $\CO_{-\kappa}$ induces the one on $\tCO_{-\kappa}$, and one can construct
an equivalence
$$
\tCO_\kappa\iso\tCO_{-\kappa}
$$
respecting bbrt structure (see ~\cite{f2}).
D.Kazhdan and G.Lusztig have constructed an equivalence
$$
\CO_{-\kappa}\iso\ _R\CC_\zeta
$$
(notations of ~\ref{br}) respecting bbrt structure (see
{}~\cite{kl} and ~\cite{l3}).
Here $\zeta=\exp(\frac{\pi\sqrt{-1}}{d\kappa})$ where $d=\max_{i\in I}d_i$.
Thus $l=2d\kappa$, and $\ell=d\kappa$.

Note that the alcoves $\nabla_\kappa$ and $\Delta_l$ (see ~\ref{alcove})
coincide.

The Kazhdan-Lusztig equivalence induces an equivalence
$$
\tCO_{-\kappa}\iso\tCO_{\zeta}
$$
where $\tCO_{\zeta}$ is the semisimple subcategory of $_R\CC_\zeta$
formed by direct sums of simple $U_B$-modules $\tL(\lambda)$
with $\lambda\in\Delta$ (see ~\cite{an} and ~\cite{ap}).
 The bbrt structure on $_R\CC_\zeta$
induces the one on $\tCO_{\zeta}$, and the last equivalence respects bbrt
structure. We denote the composition
of the above equivalences by
$$
\phi:\tCO_\kappa\iso\tCO_{\zeta}.
$$

Given any bbrt category $\CB$ and objects $L_1,\ldots, L_m\in\CB$ we obtain
a local system $\Hom_{\CB}(\One,L_1\otimes\ldots\otimes L_m)^{\heartsuit}$ on
$T\CP^{om}$ with monodromies induced by the action of braiding and balance on
$L_1\otimes\ldots\otimes L_m$.

Here and below we write
a superscript $X^{\heartsuit}$ to denote a local system over $T\CP^{om}$
with the fiber at a standard real point
$z_1<\ldots<z_m$ with tangent vectors looking to the right, equal to $X$.

Thus, given $L_1,\ldots,L_m\in\tCO_\kappa$, the local system
$$
\Hom_{\tCO_\kappa}(\One,L_1\totimes\ldots\totimes L_m)^{\heartsuit}
$$
called {\em local system of conformal blocks} is isomorphic
to the local system
$\Hom_{\tCO_{\zeta}}(\One,\phi(L_1)\totimes\ldots\totimes \phi(L_m))
^{\heartsuit}$.
Here $\totimes$ will denote the tensor product in "tilded" categories.

To unburden the notations we leave out the subscript $\zeta$ in
$_R\CC_\zeta$ from now on.

For an object $X\in\ _R\CC$ let us define  a vector space
$\langle X\rangle_{_R\CC}$ in the same manner as in ~\ref{triv sub}, i.e.
as an image of the canonical map from the maximal trivial
subobject of $X$ to the maximal trivial quotient of $X$. Given $X_1,\ldots,
X_m\in\ _R\CC$, we denote
$$
\langle X_1,\ldots X_m\rangle :=
\langle X_1\otimes\ldots\otimes X_m\rangle_{_R\CC}.
$$

\subsubsection{} {\bf Lemma.} {\em We have an isomorphism of local systems
$$
\Hom_{\tCO_{\zeta}}(\One,\phi(L_1)\totimes\ldots\totimes\phi(L_m))
^{\heartsuit}\cong
\langle\phi(L_1),\ldots,\phi(L_m)\rangle_{_R\CC}^{\heartsuit}
$$}

{\bf Proof.} Follows from ~\cite{an}. $\Box$

\subsection{Lemma}
\label{compar conf} {\em The restriction functor
$\Upsilon:\ _R\CC\lra\CC$ (cf. ~\ref{ups}) induces isomorphism
$$
\langle\phi(L_1),\ldots,\phi(L_m)\rangle_{_R\CC}\iso
\langle\Upsilon\phi(L_1),\ldots,\Upsilon\phi(L_m)\rangle_{\CC}.
$$}

{\bf Proof.} We must prove that if $\lambda_1,\ldots,\lambda_m\in\Delta$,
$\wt{L(\lambda_1)},\ldots,\wt{L(\lambda_m)}$ are corresponding simples
in $_R\CC$, and $L(\lambda_i)=\Upsilon\wt{L(\lambda_i)}$ ---
the corresponding simples in $\CC$, then
the maximal trivial direct summand of
$\wt{L(\lambda_1)}\otimes\ldots\otimes\wt{L(\lambda_m)}$ in $_R\CC$
maps isomorphically to the maximal trivial direct summand of
$L(\lambda_1)\otimes\ldots\otimes L(\lambda_m)$ in $\CC$.

According to ~\cite{an}, ~\cite{ap},
$\wt{L(\lambda_1)}\otimes\ldots\otimes\wt{L(\lambda_m)}$ is a direct sum
of a module
$\wt{L(\lambda_1)}\totimes\ldots\totimes\wt{L(\lambda_m)}\in\CO_{\zeta}$ and
a negligible module $N\in\ _R\CC$. Here {\em negligible} means
that any endomorphism of $N$ has quantum trace zero (see {\em loc. cit.}).
Moreover, it is proven in {\em loc cit.} that $N$ is a direct summand
of $W\otimes M$ for some $M\in\ _R\CC$ where
$W=\oplus_{\omega\in\Omega}\ \wt{L(\omega)}$,
$$
\Omega=\{\omega\in X|\ \langle i,\omega+\rho\rangle>0\
\mbox{for all}\ i\in I;\ \langle\beta_0,\omega+\rho\rangle=\kappa\}
$$
being the affine wall of the first alcove.
By {\em loc. cit.}, $W$ is negligible. Since
$\Upsilon\wt{L(\omega)}=L(\omega)$, $\omega\in\Omega$ and since
$\Upsilon$ commutes with braiding, balance and rigidity, we see that
the modules $L(\omega)$ are negligible in $\CC$. Hence $\Upsilon W$ is
negligible, and $\Upsilon W\otimes\Upsilon M$ is negligible, and
finally $\Upsilon N$ is negligible. This implies that $\Upsilon N$ cannot have
trivial summands (since $L(0)$) is not negligible).

We conclude that
\begin{eqnarray}\nonumber
\langle\Upsilon(\wt{L(\lambda_1)}\otimes\ldots\otimes \wt{L(\lambda_m)})\rangle
_{\CC} =
\langle\Upsilon\wt{L(\lambda_1)}\totimes\ldots\totimes
\Upsilon\wt{L(\lambda_m)}\rangle_{\CC}=\\ \nonumber
\langle\wt{L(\lambda_1)}\totimes\ldots\totimes\wt{L(\lambda_m)}\rangle_{_R\CC}
=
\langle\wt{L(\lambda_1)}\otimes\ldots\otimes\wt{L(\lambda_m)}\rangle_{_R\CC}
\ \Box\nonumber
\end{eqnarray}

\subsection{} Corollary ~\ref{conf subq} implies that the local system
$$
\langle\Upsilon\phi(L_1),\ldots,\Upsilon\phi(L_m)\rangle_{\CC}^{\heartsuit}
$$
is canonically a subquotient of the local system
$$
\Tor^{\CC}_{\infh+0}(B,\Upsilon\phi(L_1)\otimes\ldots\otimes
\Upsilon\phi(L_m)\otimes L(2\rho_\ell)^{\heartsuit}
$$
(the action of monodromy being induced by braiding
and balance on the first $m$ factors).

\subsection{} Let us fix a point $\infty\in\BP^1$ and a nonzero tangent
vector $v\in T_{\infty}\BP^1$.
This defines an open subset
$$
T\CA^{om}\subset T\CP^{om}
$$
and the locally closed embedding
$$
\xi: T\CA^{om}\hra T\CP^{om+1}.
$$
Given $\lambda_1,\ldots,\lambda_m\in\Delta$, we consider the integrable
$\hfg$-modules $\hL(\lambda_1),\ldots,\hL(\lambda_m)$ of central charge
$\kappa-\check{h}$.

Suppose that
$$
\lambda_1+\ldots +\lambda_m=\alpha\in\BN[I]\subset X.
$$
We define $\lambda_{\infty}:=2\rho_\ell$, and
$\vlambda:=(\lambda_1,\ldots,\lambda_m,\lambda_{\infty})$.
Note that $\vlambda$ is positive and
$\alpha=\alpha(\vlambda)$, in the notations of ~\ref{pos}.

Denote by $\CX^{\alpha}_{\vlambda}$ the sheaf on $\CP_{\vlambda}^{\alpha}$
obtained by gluing $\CL(\lambda_1),\ldots,\CL(\lambda_m),
\CL(\lambda_{\infty})$. Note that
$$
\CX^{\alpha}_{\vlambda}=j_{!*}\CI^{\alpha}_{\vlambda}
$$
where $j:\CP^{o\alpha}_{\vlambda}\hra \CP^{\alpha}_{\vlambda}$.

Consider the local system of conformal blocks
$$
\Hom_{\tCO_\kappa}(\One,\hL(\lambda_1)\totimes\ldots\totimes\hL(\lambda_m))
^{\heartsuit}.
$$
If $\sum_{i=1}^m\lambda_i\not\in\BN[I]\subset X$ then it vanishes by the above
comparison with its "quantum group" incarnation.

\subsection{Theorem} {\em Suppose that $\sum_{i=1}^m\lambda_i=\alpha\in\BN[I]$.
Then the local system of conformal blocks restricted to
$T\CA^{om}$ is isomorphic to a canonical subquotient of a
"geometric" local system
$$
\xi^*R^{-2m-2}\breta^{\alpha}_{m+1*}j_{!*}\CI^{\alpha}_{\vlambda}.
$$}

{\bf Proof.} This follows from Theorem ~\ref{global thm} and the previous
discussion. $\Box$

\subsection{Corollary} {\em The above local system of conformal blocks
is semisimple. It is a direct summand of the geometric local system
above.}

{\bf Proof.} The geometric system is semisimple by Decomposition theorem,
{}~\cite{bbd}, Th\'{e}or\`{e}me 6.2.5. $\Box$

\subsection{} Example ~\ref{ex conf} shows that in general
a local system of conformal blocks is a {\em proper} direct summand
of the corresponding geometric system.


\begin{thebibliography}{MMMMM}


\bibitem[A]{an} H.Andersen, Representations of quantum groups,
invariants of $3$-manifolds and semisimple tensor categories,
{\em Isr. Math. Conf. Proc.,} {\bf 7}(1993), 1-12.

\bibitem[AJS]{ajs} H.Andersen, J.Jantzen, W.Soergel, Representations of
quantum groups at $p$-th root of unity and of semisimple groups in
characteristic $p$: independence of $p$, {\em Ast\'{e}risque} {\bf 220}(1994).


\bibitem[AP]{ap} H.Andersen, J.Paradowski, Fusion categories arising
from semisimple Lie algebras, {\em CMP} {\bf 169} (1995), 563-588.


\bibitem[Ar]{a} S.Arkhipov, Semiinfinite cohomology of quantum groups,
Preprint (1995), Moscow.

\bibitem[B]{b} A.Beilinson, On the derived category of perverse sheaves,
in: K-theory, Arithmetic and Geometry, Yu.I.Manin (Ed.),
{\em Lect. Notes in Math.,} {\bf 1289}(1987), 27-41.


\bibitem[BBD]{bbd} A.Beilinson, J.Bernstein, P.Deligne, Faisceaux Pervers,
{\em Ast\'{e}risque} {\bf 100}(1982).


\bibitem[D1]{d1} P.Deligne, Une description de cat\'{e}gorie tress\'{e}e
(inspir\'{e} par Drinfeld), Letter to V.Drinfeld (1990).
Lectures at IAS (1990), unpublished.





\bibitem[F]{f2} M.Finkelberg, An equivalence of fusion categories,
{\em J. Geom. and Funct. Anal.,} to appear.

\bibitem[FS]{fs} M.Finkelberg, V.Schechtman, Localization of $\fu$-modules. I.
Intersection cohomology of real arrangements, Preprint hep-th/9411050
(1994), 1-23; II. Configuration spaces and quantum groups,
Preprint q-alg/9412017 (1994), 1-59; III. Tensor categories arising
from configuration spaces, Preprint q-alg/9503013 (1995), 1-59.

\bibitem[KL]{kl} D.Kazhdan, G.Lusztig, Tensor structures arising from
affine Lie algebras. I-IV, {\em Amer. J. Math.}, {\bf 6}(1993), 905-947;
{\bf 6}(1993), 949-1011; {\bf 7}(1994), 335-381; {\bf 7}(1994), 383-453.


\bibitem[L1]{l} G.Lusztig, Introduction to quantum groups, Boston, Birkhauser,
1993.


\bibitem[L2]{l2} G.Lusztig, Quantum groups at roots of 1, {\em Geom.
Dedicata} {\bf 35}(1990), 89-114.

\bibitem[L3]{l3} G.Lusztig, Monodromic systems on affine flag manifolds,
{\em Proc. R. Soc. Lond. A} {\bf 445} (1994), 231-246.

\bibitem[MS]{ms} G.Moore, N.Seiberg, Classical and quantum conformal field
theory, {\em Comm. Math. Phys.}, {\bf 123} (1989), 177-254.





\bibitem[Xi]{x} Xi Nanhua, Representations of finite dimensional
Hopf algebras arising from quantum groups, Preprint (1989).

\end{thebibliography}
\end{document}